\documentclass{amsart}

\newtheorem{theorem}{Theorem}[section]
\newtheorem{lemma}[theorem]{Lemma}
\newtheorem{proposition}[theorem]{Proposition}

\theoremstyle{definition}

\newtheorem{remark}[theorem]{Remark}

\usepackage{amscd,amssymb}
\def\boxtext#1{%
\vbox{%
\hrule
\hbox{\strut \vrule{} #1 \vrule}%
\hrule
}%
}

\begin{document}

\title[Generators of Invariants of Two $4 \times 4$ Matrices]
{Generators of Invariants of Two $4 \times 4$ Matrices}

\author[Vesselin Drensky and Liliya Sadikova]
{Vesselin Drensky and Liliya Sadikova}
\address{Institute of Mathematics and Informatics,
Bulgarian Academy of Sciences,
          1113 Sofia, Bulgaria}
\email{drensky@math.bas.bg}
\address{Fachgruppe Informatik,
RWTH Aachen, 52056 Aachen, Germany}
\email{sadikova@stce.rwth-aachen.de}

\thanks{The first author was partially supported by Grant MM1106/2001
of the Bulgarian National Science Fund.}
\subjclass[2000]
{16R30}
\keywords{Matrix invariants, trace algebras.}

\begin{abstract}
Over a field of characteristic 0, the algebra of invariants of
several $n\times n$ matrices under simultaneous conjugation by
$GL_n$ is generated by traces of products of generic matrices. In
this paper we have found, in terms of representation theory of
$GL_2$, a minimal set of generators of the algebra of invariants
of two $4\times 4$ matrices. The proof is purely combinatorial and
involves computer calculations with standard functions of Maple.
\end{abstract}

\maketitle

\section*{Introduction}

Let $K$ be a field of characteristic 0 and let
$C_{nd}$ be the pure trace algebra
generated by all traces of products
$\text{\rm tr}(X_{i_1}\cdots X_{i_k})$ of the generic
$n\times n$ matrices $X_i=\left(x_{pq}^{(i)}\right)$,
$p,q=1,\ldots,n$, $i=1,\ldots,d$.
The algebra $C_{nd}$ coincides with the algebra of invariants of the general
linear group $GL_n=GL_n(K)$ acting by simultaneous conjugation on
$d$ matrices of size $n\times n$.

For a background on the algebras of matrix invariants
see, e.g. \cite{F, DF}.
Traditionally, a result giving the explicit generators of the
algebra of invariants of a linear group $G$ is called a
first fundamental theorem of the invariant theory of $G$
and a result describing the relations between the generators a
second fundamental theorem. Classical
invariant theory implies that $C_{nd}$ is finitely generated.
By the Noether normalization theorem $C_{nd}$ contains
a homogeneous set of algebraically independent elements
$\{f_1,\ldots,f_k\}$ called a homogeneous system of
parameters, where $k=(d-1)n^2+1$ is the transcendence
degree of the quotient field of $C_{nd}$, such that
$C_{nd}$ is integral over the polynomial algebra $S=K[f_1,\ldots,f_k]$.
A more precise result, see Van den Bergh \cite{V},
gives that $C_{nd}$ is a finitely generated free module over
a suitable chosen $S$.
An upper bound for the generating set of the algebra $C_{nd}$
is given in terms of PI-algebras.
By the Nagata-Higman theorem the polynomial
identity $x^n=0$ implies the identity $x_1\cdots x_N=0$ for some
$N=N(n)$. If $N$ is minimal with this property, then $C_{nd}$ is
generated by traces of products $\text{\rm tr}(X_{i_1}\cdots X_{i_k})$
of degree $k\leq N$. This estimate is sharp if $d$ is sufficiently
large. The values of $N(n)$ are between $n(n+1)/2$ and $n^2$, bounds obtained
respectively by Kuzmin \cite{Ku} and Razmyslov \cite{R} (see also
\cite{DF} for an exposition of the results in \cite{Ku}).
The only explicitly known values of $N(n)$ are
$N(2)=3$ and $N(3)=6$ (an old result of Dubnov \cite{Du}), and $N(4)=10$
(a result of Vaughan-Lee \cite{VL}). Recently, Shestakov and Zhukavets
\cite{SZ} obtained that the class of nilpotency of 2-generated nil algebras
of nil index 5 is 15. All these results agree with the conjecture of Kuzmin
\cite{Ku} that $N(n)=n(n+1)/2$.

Explicit minimal sets of generators of $C_{nd}$ are found in few cases only.
By a theorem of Sibirskii \cite{S}, $C_{2d}$ is generated
by $\text{\rm tr}(X_i)$, $1\leq i\leq d$, $\text{\rm tr}(X_iX_j)$,
$1\leq i\leq j\leq d$,
$\text{\rm tr}(X_iX_jX_k)$, $1\leq i<j<k\leq d$.
Teranishi \cite{T1} found the following system of generators of $C_{32}$:
\[
\begin{array}{c}
\text{\rm tr}(X),\text{\rm tr}(Y),\text{\rm tr}(X^2),\text{\rm tr}(XY),\text{\rm tr}(Y^2), \\
\\
\text{\rm tr}(X^3),\text{\rm tr}(X^2Y),\text{\rm tr}(XY^2),\text{\rm tr}(Y^3),
\text{\rm tr}(X^2Y^2),\text{\rm tr}(X^2Y^2XY),
\end{array}
\]
where $X$, $Y$ are generic $3 \times 3$ matrices. He showed that
the first ten of these generators form a homogeneous system of
parameters of $C_{32}$ and $C_{32}$ is a free module with
generators 1 and $\text{\rm tr}(X^2Y^2XY)$ over the polynomial
algebra on these ten elements.
Abeasis and Pittaluga \cite{AP} found
a system of generators of $C_{3d}$, for any $d\geq 2$, in terms of
representation theory of the symmetric and general linear groups,
in the spirit of its use in theory of PI-algebras.
Teranishi \cite{T1, T2} found also
a homogeneous system of parameters of $C_{42}$
and a set of generators, all of them being traces of
products of generic matrices.

It has turned out that the systems of generators of Sibirskii \cite{S}
for $C_{2d}$ and of Teranishi \cite{T1} for $C_{32}$ are not very convenient
if we want to find the defining relations of the algebra.
(Compare the relation for $C_{23}$ found by Sibirskii \cite{S}
with the defining relations for $C_{2d}$ of one of the authors \cite{D2}
with respect to another natural system of generators.)
Also, it follows from the description of the generators of $C_{32}$, that
$\text{\rm tr}(X^2Y^2XY)$ satisfies a quadratic equation with coefficients
depending on the other ten generators. The explicit (but very complicated)
form of the equation was found by Nakamoto \cite{N}, over $\mathbb Z$,
with respect to a slightly different system of generators.
A much simple relation with respect to the system of generators
of $C_{32}$ in \cite{AP}
was found by Aslaksen, Drensky and Sadikova in \cite{ADS}.

The purpose of the present paper is to find a minimal system of generators
of $C_{42}$ in the spirit of the work of Abeasis and Pittaluga \cite{AP}
used in their description of the generators of $C_{nd}$.
Applied to $C_{42}$, the idea in \cite{AP} is the following.
The general linear group $GL_2$ acts
canonically on the free associative algebra $K\langle x,y\rangle$
and this induces an action on the pure trace algebra
\[
P(x,y)=K[\text{tr}(z_1\cdots z_k)\mid z_i=x,y,\ i=1,\ldots,k,\ k=1,2,\ldots],
\]
and on the algebra $C_{42}$ (because $C_{42}$
is a homomorphic image of $P(x,y)$ under the natural
homomorphism extending the map $x\to X$, $y\to Y$, where $X$ and $Y$
are two generic $4\times 4$ matrices).
The algebra $C_{42}$ has a system of generators of degree $\leq 10$. Without
loss of generality we may assume that this system consists of traces of products
$\text{tr}(Z_1\cdots Z_k)$, $Z_i=X,Y$.
Let $T_k$ be the subalgebra of $C_{42}$ generated by all traces
$\text{\rm tr}(Z_1\cdots Z_l)$ of degree $l\leq k$, $Z_i=X,Y$.
Clearly, $T_k$ is also a $GL_2$-submodule of $C_{42}$.
Let $C_{42}^{(k+1)}$ be the homogeneous component of degree $k+1$
of $C_{42}$. Then the intersection $T_k\cap C_{42}^{(k+1)}$
is a $GL_2$-module and has a complement
$G_{k+1}$ in $C_{42}^{(k+1)}$,
which is the $GL_2$-module of the ``new'' generators of degree $k+1$.
We may assume that $G_{k+1}$ is a submodule of the $GL_2$-module
spanned by traces of products $\text{\rm tr}(Z_1\cdots Z_{k+1})$ of degree $k+1$.
The $GL_2$-module of the generators of $C_{42}$ is
\[
G=G_1\oplus G_2\oplus \cdots\oplus G_{10}.
\]
The authors of \cite{AP} work with the multilinear elements
in the free algebra which are consequences of the
polynomial identity $x^n=0$. Instead, we prefer to work directly
in the algebras $P(x,y)$ and $C_{42}$ with more essential use of representation
theory of the general linear group, in the spirit of \cite{D1}.
Our main result is that the minimal generating $GL_2$-module $G$ of $C_{42}$ decomposes as
\[
\begin{array}{c}
G=W(1,0)\oplus W(2,0)\oplus W(3,0)\oplus W(4,0)\oplus W(2,2)\oplus W(3,2)\\
\oplus W(4,2)\oplus W(3,3)\oplus W(4,3)\oplus W(5,3)\oplus W(4,4)\oplus W(6,3)
\oplus W(5,5),\\
\end{array}
\]
where $W(\lambda_1,\lambda_2)$ is the irreducible $GL_2$-module
corresponding to the partition $(\lambda_1,\lambda_2)$. For
$\lambda=(\lambda_1,\lambda_2)\not=(5,5)$, one may choose as a generator of
$W(\lambda_1,\lambda_2)$ the element
\[
w_{\lambda}(X,Y)=\text{\rm tr}((XY-YX)^{\lambda_2}X^{\lambda_1-\lambda_2}).
\]
For $\lambda=(5,5)$ a generator of $W(5,5)$ is
\[
w_{(5,5)}(X,Y)=\text{\rm tr}((XY-YX)^3(X^2Y^2-XYYX-YXXY+Y^2X^2)).
\]
It seems to us that our generators are more convenient than the usual ones,
if we want to find a system of defining relations of $T_{42}$.
Most of the computations are performed with Maple, using standard functions only.
Some of the obtained relations and technical results may be used for concrete
calculations with traces of generic $4\times 4$ matrices.

\section{Preliminaries}

The group $GL_2$ acts in a canonical way on the vector space with
basis $\{x,y\}$ and this action induces a diagonal action on
the free algebra $K\langle x,y\rangle$:
\[
g(z_1\cdots z_k)=g(z_1)\cdots g(z_k),\quad z_i=x,y,\quad g\in GL_2.
\]
The action of $GL_2$ on $K\langle x,y\rangle$ induces an action on
the pure trace algebra $P(x,y)$:
\[
g(\text{tr}(z_1\cdots z_k))=\text{tr}(g(z_1)\cdots g(z_k)),\quad
z_i=x,y,\quad g\in GL_2.
\]
If we fix the generic $n\times n$ matrices $X,Y$,
the algebra of invariants $C_{n2}$ becomes a homomorphic image
of $P(x,y)$ under the natural homomorphism extending the map $x\to X$,
$y\to Y$ and the $GL_2$ action is transferred on $C_{n2}$.
The $GL_2$-module $K\langle x,y\rangle$ is a direct sum of
irreducible polynomial modules, described
in terms of partitions $\lambda=(\lambda_1,\lambda_2)$. We denote
by $W(\lambda)$ the corresponding $GL_2$-module.

The $GL_2$-submodules and factor modules $W$ of $K\langle x,y\rangle$
inherit its natural bigrading which counts the entries of $x$ and $y$
in each monomial. We denote by $W^{(p,q)}$ the corresponding homogeneous component.
The formal power series
\[
H(W,t,u)=\sum_{p,q\geq 0}\text{\rm dim}(W^{(p,q)})t^pu^q
\]
is called the Hilbert series of $W$. The Hilbert series of $W(\lambda)$ is the Schur
function $S_{\lambda}(t,u)$ which, in the case of two variables, has the very simple form
\[
S_{\lambda}(t,u)=(tu)^{\lambda_2}(t^{\lambda_1-\lambda_2}+t^{\lambda_1-\lambda_2-1}u+
\cdots+tu^{\lambda_1-\lambda_2-1}+u^{\lambda_1-\lambda_2}).
\]
The Hilbert series of $W$ plays the role of its character, and
\[
W=\bigoplus W^{\oplus m(\lambda)}(\lambda),\quad
m(\lambda)\in {\mathbb N}\cup\{0\},
\]
i.e. $W(\lambda)$ participates in $W$ with multiplicity $m(\lambda)$, if and only if
\[
H(W,t,u)=\sum m(\lambda)S_{\lambda}(t,u).
\]
The Hilbert series of $C_{42}$ was calculated by Teranishi \cite{T2},
with some misprints in the expression, and, corrected, by Berele and Stembridge \cite{BS}:
\begin{equation}\label{Hilbert series of C42}
H(C_{42},t,u)=\frac{P_C(t,u)}{(1-t)(1-u)Q_C(t,u)},
\end{equation}
\[
P_C(t,u)=(1-e_2+e_2^2)(1-e_1e_2+e_1e_2^2+e_1^2e_2^2+e_1e_2^3-e_1e_2^4+e_2^6),
\]
\[
e_1=t+u,\quad e_2=tu,
\]
\[
Q_C(t,u)=(1-t^2)(1-t^3)(1-t^4)(1-u^2)(1-u^3)(1-u^4)
\]
\[
(1-tu)^2(1-t^2u)^2(1-tu^2)^2(1-t^3u)(1-tu^3)(1-t^2u^2).
\]
We shall need the decomposition as a sum of Schur functions
of the first 11 homogeneous components $h_n(t,u)$
of the symmetric formal power series
\begin{equation}\label{Hilbert series of C0}
h(t,u)=\frac{P_C(t,u)}{Q_C(t,u)}=(1-t)(1-u)H(C_{42},t,u)
=\sum_{n\geq 0}h_n(t,u).
\end{equation}
Using the formula
\begin{equation}\label{geometric progression}
\frac{1}{1-z}=1+z+z^2+\cdots,
\end{equation}
direct calculations show that
\begin{equation}
\begin{array}{c}\label{the first eleven coefficients}
h_0=1,\quad h_1=0,\quad h_2=t^2+tu+u^2=S_{(2,0)}(t,u),\\
h_3=S_{(3,0)}(t,u),\quad h_4=2S_{(4,0)}(t,u)+2S_{(2,2)}(t,u),\\
h_5=S_{(5,0)}(t,u)+S_{(4,1)}(t,u)+2S_{(3,2)}(t,u),\\
h_6=3S_{(6,0)}(t,u)+S_{(5,1)}(t,u)+5S_{(4,2)}(t,u)+S_{(3,3)}(t,u),\\
h_7=2S_{(7,0)}(t,u)+2S_{(6,1)}(t,u)+5S_{(5,2)}(t,u)+4S_{(4,3)}(t,u),\\
h_8=4S_{(8,0)}(t,u)+2S_{(7,1)}(t,u)+10S_{(6,2)}(t,u)+6S_{(5,3)}(t,u)+8S_{(4,4)}(t,u),\\
h_9=3S_{(9,0)}(t,u)+3S_{(8,1)}(t,u)+10S_{(7,2)}(t,u)+13S_{(6,3)}(t,u)+8S_{(5,4)}(t,u),\\
h_{10}=5S_{(10,0)}(t,u)+4S_{(9,1)}(t,u)+16S_{(8,2)}(t,u)\\
+16S_{(7,3)}(t,u)+24S_{(6,4)}(t,u)+5S_{(5,5)}(t,u).\\
\end{array}
\end{equation}

Any submodule $W(\lambda)=W(\lambda_1,\lambda_2)$ of $K\langle x,y\rangle$
is generated by a unique, up to a multiplicative constant,
homogeneous element $w_{\lambda}(x,y)$ of degree $\lambda_1$ and $\lambda_2$
with respect to $x$ and $y$, respectively,
called the highest weight vector of $W(\lambda)$. It is characterized by
the following property, see \cite{DEP, ADF} and
\cite{K} for the version which we need. We state it
for two variables only.
The homogeneous polynomial $w_{\lambda}(x,y) \in K\langle x,y\rangle$ of degree
$(\lambda_1,\lambda_2)$ is a highest weight vector for some $W(\lambda_1,\lambda_2)$
if and only if $\delta(w_{\lambda}(x,y))=0$, where $\delta$ is the derivation
of $K\langle x,y\rangle$ defined by $\delta(x)=0$, $\delta(y)=x$.
If $W_i$, $i=1,\ldots,k$,
are $k$ isomorphic copies of $W(\lambda)$ and $w_i\in W_i$ are highest weight
vectors, then the highest weight vector of any submodule $W(\lambda)$
of the direct sum $W_1\oplus\cdots\oplus W_k$ has the form
$\xi_1w_1+\cdots+\xi_kw_k$ for some $\xi_i\in K$.
Any $k$ linearly independent highest weight vectors can serve
as generators of the $GL_2$-module $W_1\oplus\cdots\oplus W_k$.
Of course, similar arguments work for submodules and homomorphic images
of the $GL_2$-module $K\langle x,y\rangle$.

The multiplicity of $W(\lambda)$ in $K\langle x,y\rangle$ is equal to the degree
$d(\lambda)$ of the corresponding irreducible representation of the symmetric group
$S_k$, $k=\lambda_1+\lambda_2$,
or, equivalently, to the number of standard $\lambda$-tableaux.
For the standard $\lambda$-tableau

\bigskip
\centerline{\vbox{\hbox{$T_{\sigma}$\hskip0.5truecm$=$\hskip0.5truecm}
\hbox{\phantom{X}}} \hskip0.5truecm \vbox{\offinterlineskip
\hbox{\boxtext{$\sigma(1)$}\boxtext{$\cdots
$}\boxtext{$\sigma(2\lambda_2
-1)$}\boxtext{$\sigma(2\lambda_2+1)$}\boxtext{$\cdots
$}\boxtext{$\sigma(k)$}}
\hbox{\boxtext{$\sigma(2)$}\boxtext{$\cdots
$}\boxtext{$\phantom{-}\sigma(2\lambda_2)\phantom{-}$}}}}
\bigskip

\noindent corresponding to $\sigma\in S_k$, $k=\lambda_1+\lambda_2$,
we associate a highest weight vector $w(T_{\sigma})$ in
$K\langle x,y\rangle$. When $\sigma=\varepsilon$
is the identity of $S_k$ we fix
\[
w(T_{\varepsilon})=(xy-yx)^{\lambda_2}x^{\lambda_1-\lambda_2}
\]
\[
=\sum_{\rho_1,\ldots,\rho_s\in S_2}
\text{\rm sign}(\rho_1\cdots\rho_s)z_{\rho_1(1)}z_{\rho_1(2)}
\cdots z_{\rho_s(1)}z_{\rho_s(2)}x^{\lambda_1-\lambda_2},
\]
where $z_1=x$, $z_2=y$, $s=\lambda_2$. For $\sigma$ arbitrary,
we define $w(T_{\sigma})$ in a similar way, but the skew-symmetries
are in positions $(\sigma(1),\sigma(2)),\ldots,(\sigma(2s-1),\sigma(2s))$
instead of the positions $(1,2),\ldots,(2s-1,2s)$ (and the positions
with fixed $x$ are $\sigma(2s+1),\ldots,\sigma(k)$ instead of
$2s+1,\ldots,k$).
For the standard $(p,q)$-tableau

\bigskip
\centerline{\vbox{\hbox{$T$\hskip0.5truecm$=$\hskip0.5truecm}
\hbox{\phantom{X}}} \hskip0.5truecm \vbox{\offinterlineskip
\hbox{\boxtext{$a_1$}\boxtext{$\cdots$}\boxtext{$a_q$}\boxtext{$a_{q+1}$}\boxtext{$\cdots$}\boxtext{$a_p$}}
\hbox{\boxtext{$b_1$}\boxtext{$\cdots$}\boxtext{$b_q$}}}}
\bigskip

\noindent we denote the corresponding highest weight vector $\text{\rm tr}(w(T))$ by
\[
w(T)=w\left(\begin{matrix}
a_1&\cdots&a_q&a_{q+1}&\cdots&a_p\\
b_1&\cdots&b_q&&\\
\end{matrix}\right).
\]

\section{The $GL_2$-modules of formal traces}

We denote by $U_n$ the vector subspace of all traces of formal
products of length $n$ in the variables $x$ and $y$
in the pure trace algebra $P(x,y)$.

\begin{proposition}
The $GL_2$-modules $U_n$ spanned by the traces of the products of
length $n\leq 10$ are decomposed as follows:
\[
U_1=W(1),\quad U_2=W(2),\quad U_3=W(3),\quad U_4=W(4)\oplus
W(2,2),
\]
\[
U_5=W(5)\oplus W(3,2),\quad U_6=W(6)\oplus 2W(4,2)\oplus W(3,3),
\]
\[
U_7=W(7)\oplus 2W(5,2)\oplus 2W(4,3),
\]
\[
U_8=W(8)\oplus 3W(6,2)\oplus 3W(5,3)\oplus 3W(4,4),
\]
\[
U_9=W(9)\oplus 3W(7,2)\oplus 6W(6,3)\oplus 4W(5,4),
\]
\[
U_{10}=W(10)\oplus 4W(8,2)\oplus 7W(7,3)\oplus 10W(6,4)\oplus
4W(5,5).
\]
\end{proposition}

\begin{proof}
We calculate the Hilbert series $H(U_n,t,u)$ of $U_n$.
Since $H(U_n,t,u)$ is a symmetric
polynomial in $t,u$, it is sufficient to count the
number $h_{nk}$ of traces of degree $k$ with respect to $x$
and $n-k$ with respect to $y$, $k\geq n-k$.
There is one trace depending on $x$ only, $\text{tr}(x^n)$,
and one trace, $\text{tr}(x^{n-1}y)$, of degree $(n-1,1)$.
Hence $h_{n0}=h_{n1}=1$ and we may consider the case $k\leq n-2$ only.
Since the traces
are invariant under cyclic permutation, we may assume that
they are of the form
\[
\text{\rm tr}(x^{a_1}y^{b_1}\cdots x^{a_p}y^{b_p}),
\]
where $a_i,b_i>0$ and $a_1+\cdots+a_p=k$, $b_1+\cdots+b_p=n-k$.
We order the words $x^{a_1}y^{b_1}\cdots x^{a_p}y^{b_p}$
comparing lexicographically the $2p$-tuples
$(a_1,\ldots,a_p,b_1,\ldots,b_p)$ and fix the representative of
$\text{\rm tr}(x^{a_1}y^{b_1}\cdots x^{a_p}y^{b_p})$ to be as high
as possible.

{\it Cases $n=1,2,3$.} These cases are obvious because the coefficients
$h_{nk}, h_{n,n-k}$ of $H(U_n,t,v)$ are equal to 1 for $k=0,1$ and
\[
H(U_1,t,v)=t+u=S_{(1)}(t,u),
\]
\[
H(U_2,t,u)=t^2+tu+u^2=S_{(2)}(t,u),
\]
\[
H(U_3,t,u)=t^3+t^2u+tu^2+u^3=S_{(3)}(t,u).
\]
This means that $U_1=W(1)$, $U_2=W(2)$ and $U_3=W(3)$.

{\it Case $n=4$.} The only elements of degree $(2,2)$
are $\text{tr}(x^2y^2), \text{tr}(xyxy)$
and this implies
\[
H(U_4,t,u)=t^4+t^3u+2t^2u^2+tu^3+u^4=S_{(4)}(t,u)+S_{(2^2)}(t,u),
\]
and hence $U_4=W(4)\oplus W(2,2)$.

{\it Case $n=5$.} The traces of degree $(3,2)$ are
$\text{tr}(x^3y^2)$, $\text{tr}(x^2yxy)$ and
\[
H(U_5,t,u)=t^5+t^4u+2t^3u^2+2t^2u^3+tu^4+u^5=S_{(5)}(t,u)+S_{(3,2)}(t,u),
\]
which gives $U_5=W(5)\oplus W(3,2)$.

{\it Case $n=6$.} The traces of degree $(4,2)$ are
$\text{tr}(x^4y^2)$, $\text{tr}(x^3yxy)$, $\text{tr}(x^2yx^2y)$,
these of degree $(3,3)$ are $\text{tr}(x^3y^3)$, $\text{tr}(x^2y^2xy)$,
$\text{tr}(x^2yxy^2)$, $\text{tr}(xyxyxy)$, and
\[
H(U_6,t,u)=t^6+t^5u+3t^4u^2+4t^3u^3+3t^2u^4+tu^5+u^6
\]
\[
=S_{(6)}(t,u)+2S_{(4,2)}(t,u)+S_{(3^2)}(t,u),
\]
\[
U_6=W(6)\oplus 2W(4,2)\oplus W(3,3).
\]

{\it Case $n=7$.} The traces of degree $(5,2)$ are
$\text{tr}(x^5y^2)$, $\text{tr}(x^4yxy)$, $\text{tr}(x^3yx^2y)$,
these of degree $(4,3)$ are $\text{tr}(x^4y^3)$, $\text{tr}(x^3y^2xy)$,
$\text{tr}(x^3yxy^2)$, $\text{tr}(x^2y^2x^2y)$, $\text{tr}(x^2yxyxy)$ and
\[
H(U_7,t,u)=t^7+t^6u+3t^5u^2+5t^4u^3+5t^3u^4+3t^2u^5+tu^6+u^7
\]
\[
=S_{(7)}(t,u)+2S_{(5,2)}(t,u)+2S_{(4,3)}(t,u),
\]
\[
U_7=W(7)\oplus 2W(5,2)\oplus 2W(4,3).
\]

{\it Case $n=8$.} The traces of degree $(6,2)$ are
\[
\text{tr}(x^6y^2),\quad \text{tr}(x^{a_1}yx^{a_2}y),\quad
(a_1,a_2)=(5,1),(4,2),(3,3),
\]
and $h_{82}=4$, the traces of degree $(5,3)$ are
\[
\text{tr}(x^5y^3),
\]
\[
\text{tr}(x^{a_1}y^{b_1}x^{a_2}y^{b_2}),\quad
(a_1,a_2)=(4,1),(3,2),\quad (b_1,b_2)=(2,1),(1,2),
\]
\[
\text{tr}(x^{a_1}yx^{a_2}yx^{a_3}y),\quad
(a_1,a_2,a_3)=(3,1,1),(2,2,1),
\]
$h_{83}=7$. Finally, the traces of degree $(4,4)$ are
\[
\text{tr}(x^4y^4), \quad
\text{tr}(x^{a_1}y^{b_1}x^{a_2}y^{b_2}),
\]
\[
(a_1,a_2)=(3,1),(b_1,b_2)=(3,1),(2,2),(1,3),
\]
\[
(a_1,a_2)=(2,2),(b_1,b_2)=(3,1),(2,2),
\]
\[
\text{tr}(x^2y^{b_1}xy^{b_2}xy^{b_3}),
\quad (b_1,b_2,b_3)=(2,1,1),(1,2,1),(1,1,2),
\]
\[
\text{tr}(xyxyxyxy).
\]
Hence $h_{84}=10$,
\[
H(U_8,t,u)=t^8+t^7u+4t^6u^2+7t^5u^3+10t^4u^4+7t^3u^5+4t^2u^6+tu^7+u^8
\]
\[
=S_{(8)}(t,u)+3S_{(6,2)}(t,u)+3S_{(5,3)}(t,u)+3S_{(4,4)}(t,u),
\]
\[
U_8=W(8)\oplus 3W(6,2)\oplus 3W(5,3)\oplus 3W(4,4).
\]

{\it Case $n=9$.} The traces of degree $(7,2)$ are
\[
\text{tr}(x^7y^2),\quad \text{tr}(x^{a_1}yx^{a_2}y),\quad
(a_1,a_2)=(6,1),(5,2),(4,3),
\]
and $h_{92}=4$, these of degree $(6,3)$ are
\[
\text{tr}(x^6y^3),
\quad
\text{tr}(x^{a_1}y^{b_1}x^{a_2}y^{b_2}),
\]
\[
(a_1,a_2)=(5,1),(4,2), (b_1,b_2)=(2,1),(1,2),
\quad (a_1,a_2)=(3,3), (b_1,b_2)=(2,1),
\]
\[
\text{tr}(x^{a_1}yx^{a_2}yx^{a_3}y),\quad
(a_1,a_2,a_3)=(4,1,1),(3,2,1),(3,1,2),(2,2,2),
\]
$h_{93}=10$. The traces of degree $(5,4)$ are
\[
\text{tr}(x^5y^4),
\quad
\text{tr}(x^{a_1}y^{b_1}x^{a_2}y^{b_2}),
\quad
(a_1,a_2)=(4,1),(3,2), (b_1,b_2)=(3,1),(2,2),(1,3),
\]
\[
\text{tr}(x^{a_1}y^{b_1}x^{a_2}y^{b_2}x^{a_3}y^{b_3}),
\]
\[
(a_1,a_2,a_3)=(3,1,1),(2,2,1),\quad
(b_1,b_2,b_3)=(2,1,1),(1,2,1),(1,1,2),
\]
\[
\text{tr}(x^2yxyxyxy),
\]
$h_{94}=14$. The Hilbert series and the decomposition of $U_9$ are,
respectively,
\[
H(U_9,t,u)=t^9+t^8u+4t^7u^2+10t^6u^3+14t^5u^4
\]
\[
+14t^4u^5+10t^3u^6+4t^2u^7+tu^8+u^9
\]
\[
=S_{(9)}(t,u)+3S_{(7,2)}(t,u)+6S_{(6,3)}(t,u)+4S_{(5,4)}(t,u),
\]
\[
U_9=W(9)\oplus 3W(7,2)\oplus 6W(6,3)\oplus 4W(5,4).
\]

{\it Case $n=10$.} The traces of degree $(8,2)$ are
\[
\text{tr}(x^8y^2),\quad \text{tr}(x^{a_1}yx^{a_2}y),\quad
(a_1,a_2)=(7,1),(6,2),(5,3),(4,4),
\]
and $h_{10,2}=5$, these of degree $(7,3)$ are
\[
\text{tr}(x^7y^3),
\quad
\text{tr}(x^{a_1}y^{b_1}x^{a_2}y^{b_2}),
\]
\[
(a_1,a_2)=(6,1),(5,2),(4,3),\quad (b_1,b_2)=(2,1),(1,2),
\]
\[
\text{tr}(x^{a_1}yx^{a_2}yx^{a_3}y),\quad
(a_1,a_2,a_3)=(5,1,1),(4,2,1),(4,1,2),(3,3,1),(3,2,2),
\]
$h_{10,3}=12$. The traces of degree $(6,4)$ are
\[
\text{tr}(x^6y^4),
\quad
\text{tr}(x^{a_1}y^{b_1}x^{a_2}y^{b_2}),
\quad
(a_1,a_2)=(5,1),(4,2), (b_1,b_2)=(3,1),(2,2),(1,3),
\]
\[
(a_1,a_2)=(3,3), (b_1,b_2)=(3,1),(2,2),
\]
\[
\text{tr}(x^{a_1}y^{b_1}x^{a_2}y^{b_2}x^{a_3}y^{b_3}),
\]
\[
(a_1,a_2,a_3)=(4,1,1),(3,2,1),(3,1,2),
(b_1,b_2,b_3)=(2,1,1),(1,2,1),(1,1,2),
\]
\[
(a_1,a_2,a_3)=(2,2,2),(b_1,b_2,b_3)=(2,1,1),
\]
\[
\text{tr}(x^{a_1}yx^{a_2}yx^{a_3}yx^{a_4}y),
\quad (a_1,a_2,a_3,a_4)=(3,1,1,1), (2,2,1,1),(2,1,2,1),
\]
$h_{10,4}=22$. Finally, the traces of degree $(5,5)$ are
\[
\text{tr}(x^5y^5),
\quad
\text{tr}(x^{a_1}y^{b_1}x^{a_2}y^{b_2}),
\]
\[
(a_1,a_2)=(4,1),(3,2), \quad
(b_1,b_2)=(4,1),(3,2),(2,3),(1,4),
\]
\[
\text{tr}(x^{a_1}y^{b_1}x^{a_2}y^{b_2}x^{a_3}y^{b_3}),
\]
\[
(a_1,a_2,a_3)=(3,1,1),(2,2,1),
\]
\[
(b_1,b_2,b_3)=(3,1,1),(1,3,1),(1,1,3),(2,2,1),(2,1,2),(1,2,2),
\]
\[
\text{tr}(x^2y^{b_1}xy^{b_2}xy^{b_3}xy^{b_4}),
\]
\[
(b_1,b_2,b_3,b_4)=(2,1,1,1),(1,2,1,1),(1,1,2,1),(1,1,1,2),
\]
\[
\text{tr}(xyxyxyxyxy),
\]
hence $h_{10,5}=26$.
The Hilbert series and the decomposition of $U_{10}$ are,
respectively,
\[
H(U_{10},t,u)=t^{10}+t^9u+5t^8u^2+12t^7u^3+22t^6u^4
\]
\[
+26t^5u^5+22t^4u^6+12t^3u^7+5t^2u^8+tu^9+u^{10},
\]
\[
=S_{(10)}(t,u)+4S_{(8,2)}(t,u)+7S_{(7,3)}(t,u)+10S_{(6,4)}(t,u)+4S_{(5,5)}(t,u),
\]
\[
U_{10}=W(10)\oplus 4W(8,2)\oplus 7W(7,3)\oplus 10W(6,4)\oplus
4W(5,5).
\]
\end{proof}

\begin{lemma}
For $\lambda=(\lambda_1,\lambda_2)\vdash n$, $n\leq 10$,
the following trace polynomials form a basis of the subspace
of $U_n$ consisting of all $\lambda$-highest weight vectors:

\noindent $\lambda=(n)$, $n\geq 1$:
\[
w=w\left(\begin{matrix}
1&\cdots&n\\
\end{matrix}\right)=\text{\rm tr}(x^n);
\]

\noindent $\lambda=(2,2)$:
\[
w=\frac{1}{2}w\left(\begin{matrix}
1&3\\
2&4\\
\end{matrix}\right)
=\frac{1}{2}\text{\rm tr}([x,y]^2)
=-\text{\rm tr}(x^2y^2)+\text{\rm tr}(xyxy);
\]

\noindent $\lambda=(3,2)$:
\[
w=w\left(\begin{matrix}
1&3&5\\
2&4&\\
\end{matrix}\right)
=\text{\rm tr}([x,y]^2x)
=-\text{\rm tr}(x^3y^2)+\text{\rm tr}(x^2yxy);
\]

\noindent $\lambda=(4,2)$:
\[
w_1=w\left(\begin{matrix}
1&3&5&6\\
2&4&&\\
\end{matrix}\right)
=\text{\rm tr}([x,y]^2x^2)=
-\text{\rm tr}(x^4y^2)+2\text{\rm tr}(x^3yxy)
-\text{\rm tr}(x^2yx^2y),
\]
\[
w_2=w\left(\begin{matrix}
1&2&5&6\\
3&4&&\\
\end{matrix}\right)=
\text{\rm tr}(x^4y^2)-\text{\rm tr}(x^2yx^2y);
\]

\noindent $\lambda=(3,3)$:
\[
w=\frac{1}{3}w\left(\begin{matrix}
1&3&5\\
2&4&6\\
\end{matrix}\right)=
\text{\rm tr}(x^2y^2xy)-\text{\rm tr}(x^2yxy^2);
\]

\noindent $\lambda=(5,2)$:
\[
w_1=w\left(\begin{matrix}
1&3&5&6&7\\
2&4&&&\\
\end{matrix}\right)
=\text{\rm tr}([x,y]^2x^3)
=-\text{\rm tr}(x^5y^2)+2\text{\rm tr}(x^4yxy)
-\text{\rm tr}(x^3yx^2y)),
\]
\[
w_2=w\left(\begin{matrix}
1&2&5&6&7\\
3&4&&&\\
\end{matrix}\right)=
\text{\rm tr}(x^5y^2)-\text{\rm tr}(x^3yx^2y);
\]

\noindent $\lambda=(4,3)$:
\[
w_1=w\left(\begin{matrix}
1&3&5&7\\
2&4&6&\\
\end{matrix}\right)=
\text{\rm tr}(x^3y^2xy)-\text{\rm tr}(x^3yxy^2),
\]
\[
w_2=w\left(\begin{matrix}
1&2&5&7\\
3&4&6&\\
\end{matrix}\right)=
-\text{\rm tr}(x^4y^3)+2\text{\rm tr}(x^3y^2xy)
\]
\[
+\text{\rm tr}(x^3yxy^2)-\text{\rm tr}(x^2y^2x^2y)
-\text{\rm tr}(x^2yxyxy);
\]

\noindent $\lambda=(6,2)$:
\[
w_1=w\left(\begin{matrix}
1&3&5&6&7&8\\
2&4&&&&\\
\end{matrix}\right)
\]
\[
=\text{\rm tr}([x,y]^2x^4)
=-\text{\rm tr}(x^6y^2)+2\text{\rm tr}(x^5yxy)
-\text{\rm tr}(x^4yx^2y),
\]
\[
w_2=w\left(\begin{matrix}
1&2&5&6&7&8\\
3&4&&&&\\
\end{matrix}\right)=
\text{\rm tr}(x^6y^2)-\text{\rm tr}(x^4yx^2y),
\]
\[
w_3=w\left(\begin{matrix}
1&3&4&6&7&8\\
2&5&&&&\\
\end{matrix}\right)=
-\text{\rm tr}(x^6y^2)+\text{\rm tr}(x^5yxy)
+\text{\rm tr}(x^4yx^2y)-\text{\rm tr}(x^3yx^3y);
\]

\noindent $\lambda=(5,3)$:
\[
w_1=w\left(\begin{matrix}
1&3&5&7&8\\
2&4&6&&\\
\end{matrix}\right)
\]
\[
=\text{\rm tr}(x^4y^2xy)-\text{\rm tr}(x^4yxy^2)
-\text{\rm tr}(x^3y^2x^2y)+\text{\rm tr}(x^3yx^2y^2),
\]
\[
w_2=w\left(\begin{matrix}
1&2&5&7&8\\
3&4&6&&\\
\end{matrix}\right)=
-\text{\rm tr}(x^5y^3)+2\text{\rm tr}(x^4y^2xy)
\]
\[
-2\text{\rm tr}(x^3y^2x^2y)+2\text{\rm tr}(x^3yx^2y^2)
-\text{\rm tr}(x^2yx^2yxy),
\]
\[
w_3=w\left(\begin{matrix}
1&3&5&6&8\\
2&4&7&&\\
\end{matrix}\right)=
\text{\rm tr}(x^4y^2xy)-\text{\rm tr}(x^4yxy^2);
\]

\noindent $\lambda=(4,4)$:
\[
w_1=\frac{1}{2}w\left(\begin{matrix}
1&3&5&7\\
2&4&6&8\\
\end{matrix}\right)=\frac{1}{2}\text{\rm tr}([x,y]^4)
=-2\text{\rm tr}(x^2y^2xyxy)+2\text{\rm tr}(x^2yxy^2xy)
\]
\[
-2\text{\rm tr}(x^2yxyxy^2)+\text{\rm tr}(x^2y^2x^2y^2)
+\text{\rm tr}(xyxyxyxy),
\]
\[
w_2=w\left(\begin{matrix}
1&2&5&7\\
3&4&6&8\\
\end{matrix}\right)=
\text{\rm tr}(x^3y^3xy)-2\text{\rm tr}(x^3y^2xy^2)
+\text{\rm tr}(x^3yxy^3)
\]
\[
-\text{\rm tr}(x^2y^2xyxy)
+2\text{\rm tr}(x^2yxy^2xy)-\text{\rm tr}(x^2yxyxy^2)
-2\text{\rm tr}(x^2y^3x^2y)+2\text{\rm tr}(x^2y^2x^2y^2),
\]
\[
w_3=\frac{1}{2}w\left(\begin{matrix}
1&2&5&6\\
3&4&7&8\\
\end{matrix}\right)=
\text{\rm tr}(x^4y^4)-2\text{\rm tr}(x^3y^2xy^2)+
\text{\rm tr}(x^2yxy^2xy)
\]
\[
-2\text{\rm tr}(x^2y^3x^2y)
+2\text{\rm tr}(x^2y^2x^2y^2);
\]

\noindent $\lambda=(7,2)$:
\[
w_1=w\left(\begin{matrix}
1&3&5&6&7&8&9\\
2&4&&&&&\\
\end{matrix}\right)
\]
\[
=\text{\rm tr}([x,y]^2x^5)
=-\text{\rm tr}(x^7y^2)+2\text{\rm tr}(x^6yxy)
-\text{\rm tr}(x^5yx^2y),
\]
\[
w_2=w\left(\begin{matrix}
1&2&5&6&7&8&9\\
3&4&&&&&\\
\end{matrix}\right)=
\text{\rm tr}(x^7y^2)-\text{\rm tr}(x^5yx^2y),
\]
\[
w_3=w\left(\begin{matrix}
1&2&4&6&7&8&9\\
3&5&&&&&\\
\end{matrix}\right)=
\text{\rm tr}(x^6yxy) -\text{\rm tr}(x^4yx^3y);
\]

\noindent $\lambda=(6,3)$:
\[
w_1=w\left(\begin{matrix}
1&3&5&7&8&9\\
2&4&6&&&\\
\end{matrix}\right)
=\text{\rm tr}([x,y]^3x^3)
=\text{\rm tr}(x^5y^2xy)-\text{\rm tr}(x^5yxy^2)
\]
\[
-\text{\rm tr}(x^4y^2x^2y)+\text{\rm tr}(x^4yx^2y^2)
-\text{\rm tr}(x^3yx^2yxy)+\text{\rm tr}(x^3yxyx^2y),
\]
\[
w_2=w\left(\begin{matrix}
1&2&5&7&8&9\\
3&4&6&&&\\
\end{matrix}\right)=
-\text{\rm tr}(x^6y^3)+2\text{\rm tr}(x^5y^2xy)
\]
\[
-2\text{\rm tr}(x^4y^2x^2y)+\text{\rm tr}(x^4yx^2y^2)
+\text{\rm tr}(x^3y^2x^3y)-\text{\rm tr}(x^3yx^2yxy),
\]
\[
w_3=w\left(\begin{matrix}
1&2&4&7&8&9\\
3&5&6&&&\\
\end{matrix}\right)
\]
\[
=\text{\rm tr}(x^5yxy^2) -\text{\rm tr}(x^5y^2xy)+\text{\rm tr}(x^4yx^2y^2)-\text{\rm tr}(x^4y^2x^2y),
\]
\[
w_4=w\left(\begin{matrix}
1&3&5&6&8&9\\
2&4&7&&&\\
\end{matrix}\right)=
\text{\rm tr}(x^5y^2xy)-\text{\rm tr}(x^5yxy^2)
+\text{\rm tr}(x^4yx^2y^2)
\]
\[
-\text{\rm tr}(x^4yxyxy)-\text{\rm tr}(x^3y^2x^3y)
+\text{\rm tr}(x^3yx^2yxy)+\text{\rm tr}(x^3yxyx^2y)
-\text{\rm tr}(x^2yx^2yx^2y),
\]
\[
w_5=w\left(\begin{matrix}
1&3&5&6&7&9\\
2&4&8&&&\\
\end{matrix}\right)
\]
\[
=\text{\rm tr}(x^5y^2xy)-\text{\rm tr}(x^5yxy^2)
+\text{\rm tr}(x^3yx^2yxy)-\text{\rm tr}(x^3yxyx^2y),
\]
\[
w_6=w\left(\begin{matrix}
1&2&4&6&8&9\\
3&5&7&&&\\
\end{matrix}\right)=
-\text{\rm tr}(x^5y^2xy)+2\text{\rm tr}(x^4yx^2y^2)
\]
\[
+\text{\rm tr}(x^4yxyxy)-\text{\rm tr}(x^3y^2x^3y)
-\text{\rm tr}(x^3yxyx^2y);
\]

\noindent $\lambda=(5,4)$:
\[
w_1=w\left(\begin{matrix}
1&3&5&7&9\\
2&4&6&8&\\
\end{matrix}\right)=
\text{\rm tr}(x^3y^2x^2y^2)-\text{\rm tr}(x^3y^2xyxy) \]
\[
+\text{\rm tr}(x^3yxy^2xy)-\text{\rm tr}(x^3yxyxy^2)-\text{\rm
tr}(x^2y^2x^2yxy)+\text{\rm tr}(x^2yxyxyxy),
\]
\[
w_2=w\left(\begin{matrix}
1&2&5&7&9\\
3&4&6&8&\\
\end{matrix}\right)
=\text{\rm tr}(x^4y^3xy)-\text{\rm tr}(x^4y^2xy^2)
\]
\[
-\text{\rm tr}(x^3y^3x^2y)+\text{\rm tr}(x^3yxyxy^2)
+\text{\rm tr}(x^2y^2x^2yxy)-\text{\rm tr}(x^2yxyxyxy),
\]
\[
w_3=w\left(\begin{matrix}
1&2&4&7&9\\
3&5&6&8&\\
\end{matrix}\right)
=\text{\rm tr}(x^4y^3xy)-\text{\rm tr}(x^4yxy^3)
\]
\[
-\text{\rm tr}(x^3y^3x^2y)+\text{\rm tr}(x^3yx^2y^3)
-\text{\rm tr}(x^3y^2xyxy)+\text{\rm tr}(x^3yxyxy^2),
\]
\[
w_4=w\left(\begin{matrix}
1&2&3&4&9\\
5&6&7&8&\\
\end{matrix}\right)
=\text{\rm tr}(x^5y^4)-\text{\rm tr}(x^4y^3xy)
\]
\[
+\text{\rm tr}(x^4y^2xy^2)-\text{\rm tr}(x^4yxy^3)
-\text{\rm tr}(x^3yxy^2xy)+\text{\rm tr}(x^2yxyxyxy);
\]

\noindent $\lambda=(8,2)$:
\[
w_1=w\left(\begin{matrix}
1&3&5&6&7&8&9&10\\
2&4&&&&&&\\
\end{matrix}\right)=
-\text{\rm tr}(x^8y^2)+2\text{\rm tr}(x^7yxy)
-\text{\rm tr}(x^6yx^2y),
\]
\[
w_2=w\left(\begin{matrix}
1&3&4&6&7&8&9&10\\
2&5&&&&&&\\
\end{matrix}\right)\]
\[
=-\text{\rm tr}(x^8y^2)+\text{\rm tr}(x^7yxy)
+\text{\rm tr}(x^6yx^2y)-\text{\rm tr}(x^5yx^3y),
\]
\[
w_3=w\left(\begin{matrix}
1&3&4&5&7&8&9&10\\
2&6&&&&&&\\
\end{matrix}\right)
\]
\[
=-\text{\rm tr}(x^8y^2)+\text{\rm tr}(x^7yxy)
+\text{\rm tr}(x^5yx^3y)-\text{\rm tr}(x^4yx^4y),
\]
\[
w_4=w\left(\begin{matrix}
1&3&4&5&6&8&9&10\\
2&7&&&&&&\\
\end{matrix}\right)
\]
\[
=-\text{\rm tr}(x^8y^2)+\text{\rm tr}(x^7yxy)
-\text{\rm tr}(x^5yx^3y)+\text{\rm tr}(x^4yx^4y);
\]

\noindent $\lambda=(7,3)$:
\[
w_1=w\left(\begin{matrix}
1&3&5&7&8&9&10\\
2&4&6&&&&\\
\end{matrix}\right),
\]
\[
w_2=w\left(\begin{matrix}
1&3&5&6&8&9&10\\
2&4&7&&&&\\
\end{matrix}\right),
\]
\[
w_3=w\left(\begin{matrix}
1&3&5&6&7&9&10\\
2&4&8&&&&\\
\end{matrix}\right),
\]
\[
w_4=w\left(\begin{matrix}
1&3&5&6&7&8&10\\
2&4&9&&&&\\
\end{matrix}\right),
\]
\[
w_5=w\left(\begin{matrix}
1&2&5&6&7&9&10\\
3&4&8&&&&\\
\end{matrix}\right),
\]
\[
w_6=w\left(\begin{matrix}
1&3&4&7&8&9&10\\
2&5&6&&&&\\
\end{matrix}\right),
\]
\[
w_7=w\left(\begin{matrix}
1&2&3&7&8&9&10\\
4&5&6&&&&\\
\end{matrix}\right);
\]

\noindent $\lambda=(6,4)$:
\[
w_1=w\left(\begin{matrix}
1&3&5&7&9&10\\
2&4&6&8&&\\
\end{matrix}\right),
\]
\[
w_2=w\left(\begin{matrix}
1&2&3&4&5&6\\
7&8&9&10&&\\
\end{matrix}\right),
\]
\[
w_3=w\left(\begin{matrix}
1&2&3&7&9&10\\
4&5&6&8&&\\
\end{matrix}\right),
\]
\[
w_4=w\left(\begin{matrix}
1&2&5&7&9&10\\
3&4&6&8&&\\
\end{matrix}\right),
\]
\[
w_5=w\left(\begin{matrix}
1&2&5&6&9&10\\
3&4&7&8&&\\
\end{matrix}\right),
\]
\[
w_6=w\left(\begin{matrix}
1&2&5&6&7&10\\
3&4&8&9&&\\
\end{matrix}\right),
\]
\[
w_7=w\left(\begin{matrix}
1&3&4&7&8&10\\
2&5&6&9&&\\
\end{matrix}\right),
\]
\[
w_8=w\left(\begin{matrix}
1&2&4&6&8&10\\
3&5&7&9&&\\
\end{matrix}\right),
\]
\[
w_9=w\left(\begin{matrix}
1&3&4&5&8&9\\
2&6&7&10&&\\
\end{matrix}\right),
\]
\[
w_{10}=w\left(\begin{matrix}
1&3&4&7&9&10\\
2&5&6&8&&\\
\end{matrix}\right);
\]

\noindent $\lambda=(5,5)$:
\[
w_1=w\left(\begin{matrix}
1&3&5&7&9\\
2&4&6&8&10\\
\end{matrix}\right),
\]
\[
w_2=w\left(\begin{matrix}
1&3&5&7&8\\
2&4&6&9&10\\
\end{matrix}\right),
\]
\[
w_3=w\left(\begin{matrix}
1&3&5&6&7\\
2&4&8&9&10\\
\end{matrix}\right),
\]
\[
w_4=w\left(\begin{matrix}
1&2&3&4&9\\
5&6&7&8&10\\
\end{matrix}\right).
\]
\end{lemma}

\begin{proof}
For each $\lambda\vdash n$, the number of highest weight vectors
$w_i$ given above, coincides with the multiplicity of $W(\lambda)$
in $U_n$. The explicit form of $w_i$ is obtained by direct
calculations. It is sufficient to prove that the $w_i$s are
linearly independent in $U_n$. For most of the cases it is almost
obvious. In the more complicated cases the linear independence is
obtained by calculating the rank of the matrix with entries equal
to the coefficients of some of $\text{tr}(x^{a_1}y^{b_1}\cdots
x^{a_p}y^{b_p})$ in the expression of $w_i$. For example, for
$\lambda=(6,3)$ we compare the coefficients of
\[
\text{\rm tr}(x^6y^3),\text{\rm tr}(x^5y^2xy),\text{\rm
tr}(x^5yxy^2),\text{\rm tr}(x^4y^2x^2y),\text{\rm
tr}(x^4yx^2y^2),\text{\rm tr}(x^3yx^2yxy)
\]
in the presentation of $w_1,\ldots,w_6$ and obtain the matrix
\[
\left(\begin{array}{rrrrrr}
0&0&1&0&1&0\\
1&1&-1&-1&-1&2\\
-1&-1&-1&0&-1&-2\\
-1&0&1&0&1&0\\
1&1&-1&2&0&0\\
-1&1&0&0&0&-1\\
\end{array}\right).
\]
Its rank is equal to 6, and this gives that the trace polynomials
$w_1,\ldots,w_6$ are linearly independent in $U_9$.
\end{proof}

\section{The new generators}

We fix the generic $4\times 4$ matrices $X=(x_{ij})$ and $Y=(y_{ij})$.
Changing the variables on the diagonals
of $X$ and $Y$, we can present them in the form
\begin{equation}\label{from generic to traceless matrices}
X=\frac{1}{4}\text{\rm tr}(X)e+x,\quad
Y=\frac{1}{4}\text{\rm tr}(Y)e+y,
\end{equation}
where $e$ is the identity $4\times 4$ matrix and $x$, $y$ are generic traceless
matrices. As in the case of generic matrices without restrictions on the trace,
we may assume that $x$ is diagonal. We fix the notation
\[
x=\left(\begin{matrix}
x_1&0&0&0\\
0&x_2&0&0\\
0&0&x_3&0\\
0&0&0&-(x_1+x_2+x_3)\\
\end{matrix}\right),
\]
\[
y=\left(\begin{matrix}
y_{11}&y_{12}&y_{13}&y_{14}\\
y_{21}&y_{22}&y_{23}&y_{24}\\
y_{31}&y_{32}&y_{33}&y_{34}\\
y_{41}&y_{42}&y_{43}&-(y_{11}+y_{22}+y_{33})\\
\end{matrix}\right),
\]
where $x_1,x_2,x_3$, $y_{ij}$ are algebraically independent
commuting variables. The algebra $C_{42}$ is generated by
$\text{tr}(X),\text{tr}(Y)$ and the traces of products of $x$ and $y$.
We denote by $C_0$ the algebra generated by $\text{tr}(z_1\cdots z_k)$,
$z_i=x,y$, and $k\geq 2$. Clearly,
\begin{equation}\label{the tensor product}
C_{42}\cong K[\text{\rm tr}(X),\text{\rm tr}(Y)]\otimes_KC_0,
\end{equation}
and it is sufficient to find a minimal system of generators of $C_0$.
Our idea is the following. We want to build inductively the
$GL_2$-module $G$ of the a minimal system of generators. We start with
the $GL_2$-module $G_1=W(1)$
generated by $\text{tr}(X)$, and further, we assume that the generators
of higher degree are in $C_0$. Assume that we
have already found the module $G_2\oplus\cdots\oplus G_{n-1}$ of generators of degree
$<n$ of $C_0$. Then we consider the symmetric algebra
$K[G_2\oplus\cdots\oplus G_{n-1}]$ and all
trace identities for $M_4(K)$ of the form $f(x,y)-g(x,y)=0$, where
$f(x,y)\in K[G_2\oplus\cdots\oplus G_{n-1}]$ and $g(x,y)\in U_n$. The elements $g(x,y)$
generate a $GL_2$-submodule $R_n$ of $U_n$. Then the
$GL_2$-submodule $G_n$ of $U_n$ is a complement of $R_n$ in $U_n$.
We may identify $G_n$ with a $GL_2$-submodule of $C_0$.
We add $G_n$ to the module of generators and
call it the $GL_2$-module of new generators of degree $n$.
Then $G_2\oplus\cdots\oplus G_{10}$ is a $GL_2$-module of
generators of the algebra of invariants of two traceless $4\times 4$
matrices and $G=G_1\oplus (G_2\oplus\cdots\oplus G_{10})$
is a minimal $GL_2$-module of generators of $C_{42}$.
Let the homogeneous component of degree $n$ of
$K[G_2\oplus\cdots\oplus G_{n-1}]$ be
\[
K[G_2\oplus\cdots\oplus G_{n-1}]^{(n)}
=\sum_{(\lambda_1,\lambda_2)\vdash n}p(\lambda)W(\lambda).
\]
Similarly, we have the decomposition
\[
U_n=\sum_{(\lambda_1,\lambda_2)\vdash n}q(\lambda)W(\lambda)
\]
from Proposition 2.1 and the decomposition
\[
R_n=\sum_{(\lambda_1,\lambda_2)\vdash n}r(\lambda)W(\lambda)\subset U_n,
\]
with unknown multiplicities $r(\lambda)$.
If $v_1,\ldots,v_{p(\lambda)}\in K[G_2\oplus\cdots\oplus G_{n-1}]^{(n)}$
are linearly independent highest weight
vectors of the direct sum of $p=p(\lambda)$ isomorphic copies of $W(\lambda)$
and $w_1,\ldots,w_q\in U_n$, $q=q(\lambda)$, are
linearly independent highest weight vectors in
$qW(\lambda)\in U_n$, then the highest weight vectors of the
$GL_2$-submodule $rW(\lambda)\subset R_n$, $r=r(\lambda)$,
are all trace identities of the form
\[
w(x,y)=(\xi_1v_1+\cdots+\xi_pv_p)-(\eta_1w_1+\cdots+\eta_qw_q)=0,
\]
where $\xi_1,\ldots,\xi_p,\eta_1,\ldots,\eta_q$ are constants from the field $K$.
We consider the elements $\xi_1,\ldots,\xi_p$, $\eta_1,\ldots,\eta_q$ as unknowns.
The entries of the evaluation of $w(x,y)$ on the generic traceless matrices $x,y$ are
equal to 0. In this way we obtain
a linear homogeneous system with respect to the unknowns
$\xi_1,\ldots,\xi_p$, $\eta_1,\ldots,\eta_q$,
and with coefficients which are polynomials in $x_1,x_2,x_3$ and
$y_{ij}$. The module $R_n$
is generated by all $\eta_1w_1+\cdots+\eta_qw_q$ which correspond
to solutions $\xi_1,\ldots,\xi_p$, $\eta_1,\ldots,\eta_q$.

The highest weight vectors $w_1,\ldots,w_q$ are given in the
previous section. In order to find $v_1,\ldots,v_p$, we need
first to know the multiplicity $p(\lambda)$ of $W(\lambda)$ in
$K[G_2\oplus\cdots\oplus G_{n-1}]^{(n)}$. If
$G_2\oplus\cdots\oplus G_{n-1}=W_1\oplus  \cdots\oplus W_m$, where the $GL_2$-modules
$W_1,\ldots,W_m$ are irreducible, then
\[
K[G_2\oplus\cdots\oplus G_{n-1}]=K[W_1]\otimes\cdots\otimes K[W_m].
\]
If $W_i$ is homogeneous of degree $d_i$, and
$W_i^{\otimes_sb_i}$ denotes the $b_i$ symmetric power of $W_i$, then
\[
K[G_2\oplus\cdots\oplus G_{n-1}]^{(n)}
=\bigoplus\left(W_1^{\otimes_sb_1}\otimes \cdots
\otimes W_m^{\otimes_sb_m}\right),
\]
where the direct sum is on all $(b_1,\ldots,b_m)$ with
$b_1d_1+\cdots+b_md_m=n$. The decomposition of $W_i^{\otimes_sb_i}$
can be obtained from the Hilbert series of $K[W_i]$ and this of the usual
tensor product by the Littlewood-Richardson rule.
From computational point of view it is convenient to
calculate the homogeneous component of degree $n$
of the Hilbert series of $K[G_2\oplus\cdots\oplus G_{n-1}]$. If
\[
H(W_i,t,u)=a_{i0}t^{d_i}+a_{i1}t^{d_i-1}u
+\cdots+a_{i,d_i-1}tu^{d_i-1}+a_{id_i}u^{d_i},\quad
a_{ij}\geq 0,
\]
is the Hilbert series of $W_i$, then
\[
H(K[G_2\oplus\cdots\oplus G_{n-1}],t,u)
=\prod_{i=1}^m\prod_{j=0}^{d_i}\frac{1}{(1-t^ju^{d_i-j})^{a_{ij}}}.
\]
Then the homogeneous component of degree $n$ of the Hilbert series
can be obtained using (\ref{geometric progression}). The highest weight
vectors can be found using the constructive approach with the derivations.
We start the realization of this scheme step by step.

\begin{lemma} For $n\leq 5$, the $GL_2$-modules $G_n$
of the new generators of $C_0$ are
\[
G_2=W(2,0),\quad G_3=W(3,0),\quad G_4=W(4,0)\oplus W(2,2),\quad
G_5=W(3,2).
\]
\end{lemma}

\begin{proof}
By Proposition 2.1
\[
U_2\oplus U_3\oplus U_4\oplus U_5=
W(2)\oplus W(3)\oplus W(4)\oplus
W(2,2)\oplus W(5,0)\oplus W(3,2).
\]
Hence $G_2\oplus G_3\oplus G_4\oplus G_5$ is a homomorphic image of this module.
We shall show that $W(5,0)$ belongs to the subalgebra generated by
$G_2\oplus G_3\oplus G_4$. Since $W(5,0)$ is generated by $\text{\rm tr}(x^5)$,
it is sufficient to show that $\text{\rm tr}(x^5)$ can be expressed in $C_0$
in terms of $\text{\rm tr}(x^2)$, $\text{\rm tr}(x^3)$, and $\text{\rm tr}(x^4)$.
By the Cayley-Hamilton theorem, $x$ satisfies an equation of degree 4
with coefficients depending on $\text{\rm tr}(x^2)$, $\text{\rm tr}(x^3)$,
and $\text{\rm tr}(x^4)$. Direct verification gives that
\[
x^4-\frac{1}{2}\text{\rm tr}(x^2)x^2-\frac{1}{3}\text{\rm tr}(x^3)x
+\left(-\frac{1}{4}\text{\rm tr}^2(x^2)+\frac{1}{8}\text{\rm tr}(x^4)\right)e=0.
\]
Multiplying by $x$, taking the trace, and using that $\text{\rm tr}(x)=0$,
we obtain
\begin{equation}\label{trace of xxxxx}
\text{\rm tr}(x^5)-\frac{5}{6}\text{\rm tr}(x^2)\text{\rm tr}(x^3)=0.
\end{equation}
(In the same way, by easy induction, we obtain that
$\text{\rm tr}(x^n)\in K[\text{\rm tr}(x^2),\text{\rm tr}(x^3),\text{\rm tr}(x^4)]$
in $C_0$ for any $n\geq 5$.) Hence,
the $GL_2$-module of the new generators of $C_0$,
of degree $\leq 5$, is a homomorphic image of the module
\[
F=W(2)\oplus W(3)\oplus W(4)\oplus
W(2,2)\oplus W(3,2).
\]
The Hilbert series of $F$ is
\[
H(F,t,u)=(t^2+tu+u^2)+(t^3+t^2u+tu^2+u^3)
\]
\[
+(t^4+t^3u+2t^2u^2+tu^3+u^4)
+(t^3u^2+t^2u^3)
\]
and the Hilbert series of the symmetric algebra of $F$ is
\[
H(K[F],t,u)=\frac{1}{q_2(t,u)q_3(t,u)q_4(t,u)q_5(t,u)},
\]
\begin{equation}
\begin{array}{c}\label{first four polynomials}
q_2(t,u)=(1-t^2)(1-tu)(1-u^2),\\
q_3(t,u)=(1-t^3)(1-t^2u)(1-tu^2)(1-u^3),\\
q_4(t,u)=(1-t^4)(1-t^3u)(1-t^2u^2)^2(1-tu^3)(1-u^4),\\
q_5(t,u)=(1-t^3u^2)(1-t^2u^3).\\
\end{array}
\end{equation}
The algebra $K[F]$ maps naturally on the subalgebra of $C_0$ generated by
$G_2\oplus G_3\oplus G_4\oplus G_5$, and the coefficients of the
Hilbert series of $K[F]$ are bigger or equal to the corresponding
coefficients of the Hilbert series of this subalgebra.
By direct calculation, we compare the coefficients of the homogeneous
components of degree $\leq 5$ of the Hilbert series of $H(K[F],t,u)$ and $H(C_0,t,u)$.
We find that they coincide and this implies that
the homogeneous components of degree $\leq 5$ of $K[F]$ and $C_0$
are isomorphic as $GL_2$-modules, completing the proof of the lemma.
\end{proof}

In the following lemmas the elements $w_i$ are taken from Lemma 2.2.
The relations have been found using the scheme described in the beginning
of the section. The proofs are direct verifications, performed by Maple.

\begin{lemma}
For $\lambda=(4,2)$, the following is a relation in $C_0$:
\[
6w_1-12w_2+6v_1+2v_2-3v_3-5v_4=0,
\]
where
\[
v_1=\text{\rm tr}(x^3)\text{\rm tr}(xy^2)-\text{\rm tr}^2(x^2y),
\]
\[
v_2=3\text{\rm tr}(x^4)\text{\rm tr}(y^2)-6\text{\rm tr}(x^3y)\text{\rm tr}(xy)
+(2\text{\rm tr}(x^2y^2)+\text{\rm tr}(xyxy))\text{\rm tr}(x^2),
\]
\[
v_3=\text{\rm tr}(x^2)(\text{\rm tr}(x^2)\text{\rm tr}(y^2)-\text{\rm tr}^2(xy)),
\]
\[
v_4=\frac{1}{2}\text{\rm tr}([x,y]^2)\text{\rm tr}(x^2).
\]
\end{lemma}

\begin{lemma}
The following relations hold in $C_0$:

For $\lambda=(5,2)$:
\[
-24w_1+18v_1+4v_2+16v_3-3v_4-4v_5=0,
\]
\[
24w_2+6v_1-4v_2+3v_4+4v_5=0,
\]
where
\[
v_1=\text{\rm tr}([x,y]^2x)\text{\rm tr}(x^2),
\]
\[
v_2=3\text{\rm tr}(x^4)\text{\rm tr}(xy^2)-6\text{\rm tr}(x^3y)\text{\rm tr}(x^2y)
+(\text{\rm tr}(xyxy)+2\text{\rm tr}(x^2y^2))\text{\rm tr}(x^3),
\]
\[
v_3=\frac{1}{2}(\text{\rm tr}([x,y]^2)\text{\rm tr}(x^3),
\]
\[
v_4=\text{\rm tr}(x^2)(\text{\rm tr}(x^2)\text{\rm tr}(xy^2)-2\text{\rm tr}(xy)\text{\rm tr}(x^2y)
+\text{\rm tr}(y^2)\text{\rm tr}(x^3),
\]
\[
v_5=(-\text{\rm tr}(x^2))\text{\rm tr}(y^2)+\text{\rm tr}^2(xy))\text{\rm tr}(x^3).
\]

For $\lambda=(4,3)$:
\[
4w_1-8w_2-10v_1-4v_2+3v_3=0,
\]
where
\[
v_1=\text{\rm tr}([x,y]^2x)\text{\rm tr}(xy)-\text{\rm tr}([x,y]^2y)\text{\rm tr}(x^2),
\]
\[
v_2=\text{\rm tr}(x^4)\text{\rm tr}(y^3)
-3\text{\rm tr}(x^3y)\text{\rm tr}(xy^2)
\]
\[
+(\text{\rm tr}(xyxy)+2\text{\rm tr}(x^2y^2))\text{\rm tr}(x^2y)
-\text{\rm tr}(xy^3)\text{\rm tr}(x^3),
\]
\[
v_3=\text{\rm tr}^2(x^2)\text{\rm tr}(y^3)-3\text{\rm tr}(x^2)\text{\rm tr}(xy)\text{\rm tr}(xy^2)
\]
\[
+\text{\rm tr}(x^2)\text{\rm tr}(y^2)\text{\rm tr}(x^2y)
+2\text{\rm tr}^2(xy)\text{\rm tr}(x^2y)-\text{\rm tr}(xy)\text{\rm tr}(y^2)\text{\rm tr}(x^3).
\]
\end{lemma}

\begin{lemma}
The following relations hold in $C_0$:

For $\lambda=(6,2)$:
\[
12w_1-6v_1-4v_2-6v_4+3v_7=0,
\]
\[
24w_2-6v_1+8v_2-6v_5-2v_6+5v_7-6v_8+6v_9+3v_{10}=0,
\]
\[
48w_3+18v_1-40v_2-16v_3-8v_4+18v_5+10v_6-13v_7+2v_8-6v_9-9v_{10}=0,
\]
where
\[
v_1=\text{\rm tr}([x,y]^2x^2)\text{\rm tr}(x^2),
\]
\[
v_2=\text{\rm tr}([x,y]^2x)\text{\rm tr}(x^3),
\]
\[
v_3=\text{\rm tr}(x^4)(2\text{\rm tr}(x^2y^2)+\text{\rm tr}(xyxy))-3\text{\rm tr}^2(x^3y),
\]
\[
v_4=\frac{1}{2}\text{\rm tr}([x,y]^2)\text{\rm tr}(x^4),
\]
\[
v_5=\text{\rm tr}(x^4)(\text{\rm tr}(x^2)\text{\rm tr}(y^2)-\text{\rm tr}^2(xy),
\]
\[
v_6=3\text{\rm tr}(x^4)\text{\rm tr}^2(xy)
-6\text{\rm tr}(x^3y)\text{\rm tr}(x^2)\text{\rm tr}(xy)
+(2\text{\rm tr}(x^2y^2)+\text{\rm tr}(xyxy))\text{\rm tr}^2(x^2),
\]
\[
v_7=\frac{1}{2}\text{\rm tr}([x,y]^2)\text{\rm tr}^2(x^2),
\]
\[
v_8=\text{\rm tr}(x^3)(\text{\rm tr}(x^3)\text{\rm tr}(y^2)
-2\text{\rm tr}(x^2y)\text{\rm tr}(xy)+\text{\rm tr}(xy^2)\text{\rm tr}(x^2)),
\]
\[
v_9=\text{\rm tr}^2(x^3)\text{\rm tr}(y^2)
-2\text{\rm tr}(x^3)\text{\rm tr}(x^2y)\text{\rm tr}(xy)
+\text{\rm tr}^2(x^2y)\text{\rm tr}(x^2),
\]
\[
v_{10}=(\text{\rm tr}(x^2)\text{\rm tr}(y^2)-\text{\rm tr}^2(xy))\text{\rm tr}^2(x^2).
\]

For $\lambda=(5,3)$:
\[
16w_1-8w_2-v_1-4v_2-4v_3+2v_4-v_5=0,
\]
\[
2w_3-v_2=0,
\]
where
\[
v_1=-2\text{\rm tr}([x,y]^2x^2)\text{\rm tr}(xy)
+(\text{\rm tr}([x,y]^2xy)+\text{\rm tr}([x,y]^2yx))\text{\rm tr}(x^2),
\]
\[
v_2=(\text{\rm tr}(x^2y^2xy)-\text{\rm tr}(y^2x^2yx))\text{\rm tr}(x^2),
\]
\[
v_3=\text{\rm tr}([x,y]^2x)\text{\rm tr}(x^2y)-\text{\rm tr}([x,y]^2y)\text{\rm tr}(x^3),
\]
\[
v_4=\text{\rm tr}(x^4)\text{\rm tr}(xy)\text{\rm tr}(y^2)
-\text{\rm tr}(x^3y)\text{\rm tr}(x^2)\text{\rm tr}(y^2)
\]
\[
-2\text{\rm tr}(x^3y)\text{\rm tr}^2(xy)
+(2\text{\rm tr}(x^2y^2)+\text{\rm tr}(xyxy)\text{\rm tr}(x^2)\text{\rm tr}(xy)
-\text{\rm tr}(xy^3)\text{\rm tr}^2(x^2),
\]
\[
v_5=-2\text{\rm tr}(x^3)\text{\rm tr}(xy^2)\text{\rm tr}(xy)
+2\text{\rm tr}x^2y)\text{\rm tr}(x^2y)\text{\rm tr}(xy)
\]
\[
+\text{\rm tr}(x^3)\text{\rm tr}(y^3)\text{\rm tr}(x^2)
-\text{\rm tr}(x^2y)\text{\rm tr}(xy^2)\text{\rm tr}(x^2).
\]

For $\lambda=(4,4)$:
\[
-24w_1+12w_2-12v_1-4v_2+16v_3+v_4+16v_5=0,
\]
\[
-72w_1+24w_3-42v_1-20v_2+56v_3+3v_4+66v_5+6v_6+6v_7=0,
\]
where
\[
v_1=\text{\rm tr}([x,y]^2x^2)\text{\rm tr}(y^2)
-(\text{\rm tr}([x,y]^2xy)+\text{\rm tr}([x,y]^2yx))\text{\rm tr}(xy)
+\text{\rm tr}([x,y]^2y^2)\text{\rm tr}(x^2),
\]
\[
v_2=3\text{\rm tr}(x^4)\text{\rm tr}(y^4)-12\text{\rm tr}(x^3y)\text{\rm tr}(xy^3)
+(2\text{\rm tr}(x^2y^2)+\text{\rm tr}(xyxy))^2,
\]
\[
v_3=\frac{1}{4}\text{\rm tr}^2([x,y]^2),
\]
\[
v_4=3\text{\rm tr}(x^4)\text{\rm tr}^2(y^2)-12\text{\rm tr}(x^3y)\text{\rm tr}(xy)\text{\rm tr}(y^2)
\]
\[
+2(2\text{\rm tr}(x^2y^2)+\text{\rm tr}(xyxy))\text{\rm tr}(x^2)\text{\rm tr}(y^2)
+4(2\text{\rm tr}(x^2y^2)+\text{\rm tr}(xyxy))\text{\rm tr}^2(xy)
\]
\[
-12\text{\rm tr}(xy^3)\text{\rm tr}(x^2)\text{\rm tr}(xy)
+3\text{\rm tr}(y^4)\text{\rm tr}^2(x^2),
\]
\[
v_5=\frac{1}{2}\text{\rm tr}([x,y]^2)(\text{\rm tr}(x^2)\text{\rm tr}(y^2)-\text{\rm tr}^2(xy)),
\]
\[
v_6=-\text{\rm tr}(x^3)\text{\rm tr}(xy^2)\text{\rm tr}(y^2)
+\text{\rm tr}^2(x^2y)\text{\rm tr}(y^2)
+\text{\rm tr}(x^3)\text{\rm tr}(y^3)\text{\rm tr}(xy)
\]
\[
-\text{\rm tr}(x^2y)\text{\rm tr}(xy^2)\text{\rm tr}(xy)
-\text{\rm tr}(x^2y)\text{\rm tr}(y^3)\text{\rm tr}(x^2)
+\text{\rm tr}^2(xy^2)\text{\rm tr}(x^2),
\]
\[
v_7=(\text{\rm tr}(x^2)\text{\rm tr}(y^2)-\text{\rm tr}^2(xy))^2.
\]
\end{lemma}

\begin{lemma} The following relations hold in $C_0$ between the
elements $w_i$ and the elements $v_j$ given below:

For $\lambda=(7,2)$:
\[
48w_1-16v_1-12v_2-12v_3-4v_5-16v_7-v_9+4v_{10}=0,
\]
\[
144w_2-24v_1+36v_2-24v_4-20v_5+48v_6+20v_7-24v_8+9v_9=0,
\]
\[
48w_3-12v_2-4v_5-8v_7-v_9+4v_{10}=0,
\]
\[
v_1=\text{\rm tr}([x,y]^2x^2)\text{\rm tr}(x^3),
\]
\[
v_2=\text{\rm tr}([x,y]^2x)\text{\rm tr}(x^4),
\]
\[
v_3=\text{\rm tr}([x,y]^2x)\text{\rm tr}^2(x^2),
\]
\[
v_4=\text{\rm tr}(x^4)(\text{\rm tr}(x^3)\text{\rm tr}(y^2)
+\text{\rm tr}(xy^2)\text{\rm tr}(x^2)-2\text{\rm tr}(x^2y)\text{\rm tr}(xy)),
\]
\[
v_5=3\text{\rm tr}(x^4)\text{\rm tr}(xy^2)\text{\rm tr}(x^2)
\]
\[
+(2\text{\rm tr}(x^2y^2)+\text{\rm tr}(xyxy))\text{\rm tr}(x^3)\text{\rm tr}(x^2)
-6\text{\rm tr}(x^3y)\text{\rm tr}(x^2y)\text{\rm tr}(x^2),
\]
\[
v_6=\text{\rm tr}(x^4)\text{\rm tr}(xy^2)\text{\rm tr}(x^2)
-\text{\rm tr}(x^4)\text{\rm tr}(x^2y)\text{\rm tr}(xy)
\]
\[
+\text{\rm tr}(x^3y)\text{\rm tr}(x^3)\text{\rm tr}(xy)
-\text{\rm tr}(x^3y)\text{\rm tr}(x^2y)\text{\rm tr}(x^2),
\]
\[
v_7=\frac{1}{2}\text{\rm tr}([x,y]^2)\text{\rm tr}(x^3)\text{\rm tr}(x^2),
\]
\[
v_8=\text{\rm tr}(x^3)(\text{\rm tr}(x^3)\text{\rm tr}(xy^2)-\text{\rm tr}^2(x^2y)),
\]
\[
v_9=(\text{\rm tr}(x^3)\text{\rm tr}(y^2)-2\text{\rm tr}(x^2y)\text{\rm tr}(xy)
+\text{\rm tr}(xy^2)\text{\rm tr}(x^2))\text{\rm tr}^2(x^2),
\]
\[
v_{10}=(\text{\rm tr}(x^3)\text{\rm tr}^2(xy)-2\text{\rm tr}(x^2y)\text{\rm tr}(x^2)\text{\rm tr}(xy)
+\text{\rm tr}(xy^2)\text{\rm tr}^2(x^2))\text{\rm tr}(x^2).
\]

For $\lambda=(6,3)$:
\[
-24w_1+48w_2-12v_1+3v_2-24v_3+36v_4+12v_5
\]
\[
-6v_6-2v_7+8v_8+7v_9+3v_{10}-15v_{11}-7v_{12}=0,
\]
\[
2w_3+v_1=0,
\]
\[
24w_4-12v_1-12v_2-12v_3-24v_4+6v_5-4v_7+4v_8-16v_9-6v_{11}-v_{12}=0,
\]
\[
w_1+w_5-v_1-v_3=0,
\]
\[
24w_1+48w_6+3v_2-24v_3-12v_4-18v_5
\]
\[
-6v_6-10v_7-4v_8-9v_9+3v_{10}+3v_{11}+6v_{12}=0,
\]
\[
v_1=\text{\rm tr}([x,y]^3x)\text{\rm tr}(x^2),
\]
\[
v_2=2\text{\rm tr}([x,y]x^2)\text{\rm tr}(x^2y)
-(\text{\rm tr}([x,y]^2xy)+\text{\rm tr}([x,y]^2yx))\text{\rm tr}(x^3),
\]
\[
v_3=(\text{\rm tr}(x^2y^2xy)-\text{\rm tr}(y^2x^2yx))(\text{\rm tr}(x^3),
\]
\[
v_4=\text{\rm tr}([x,y]^2x)\text{\rm tr}(x^3y)-\text{\rm tr}([x,y]^2y)\text{\rm tr}(x^4),
\]
\[
v_5=(\text{\rm tr}([x,y]^2x)\text{\rm tr}(xy)-\text{\rm tr}([x,y]^2y)\text{\rm tr}(x^2))\text{\rm tr}(x^2),
\]
\[
v_6=(\text{\rm tr}(x^4)\text{\rm tr}(x^2y)
-\text{\rm tr}(x^3y)\text{\rm tr}(x^3))\text{\rm tr}(y^2)
\]
\[
+2(-\text{\rm tr}(x^4)\text{\rm tr}(xy^2)
+\text{\rm tr}(x^3y)\text{\rm tr}(x^2y))\text{\rm tr}(xy)
\]
\[
+(2\text{\rm tr}(x^3y)\text{\rm tr}(xy^2)
-(2\text{\rm tr}(x^2y^2)+\text{\rm tr}(xyxy))\text{\rm tr}(x^2y)
+\text{\rm tr}(xy^3)\text{\rm tr}(x^3))\text{\rm tr}(x^2),
\]
\[
v_7=(3\text{\rm tr}(x^4)\text{\rm tr}(xy^2)
-6\text{\rm tr}(x^3y)\text{\rm tr}(x^2y)
+(2\text{\rm tr}(x^2y^2)+\text{\rm tr}(xyxy))\text{\rm tr}(x^3))\text{\rm tr}(xy)
\]
\[
+(-3\text{\rm tr}(x^3y)\text{\rm tr}(xy^2)
+2(2\text{\rm tr}(x^2y^2)+\text{\rm tr}(xyxy))\text{\rm tr}(x^2y)
-3\text{\rm tr}(xy^3)\text{\rm tr}(x^3))\text{\rm tr}(x^2),
\]
\[
v_8=(\text{\rm tr}(x^4)\text{\rm tr}(y^3)-3\text{\rm tr}(x^3y)\text{\rm tr}(xy^2)
\]
\[
+(2\text{\rm tr}(x^2y^2)+\text{\rm tr}(xyxy))\text{\rm tr}(x^2y)
-\text{\rm tr}(xy^3)\text{\rm tr}(x^3))\text{\rm tr}(x^2),
\]
\[
v_9=\frac{1}{2}\text{\rm tr}([x,y]^2)
(\text{\rm tr}(x^3)\text{\rm tr}(xy)-\text{\rm tr}(x^2y)\text{\rm tr}(x^2)),
\]
\[
v_{10}=\text{\rm tr}^2(x^3)\text{\rm tr}(y^3)
-3\text{\rm tr}(x^3)\text{\rm tr}(x^2y)\text{\rm tr}(xy^2)
+2\text{\rm tr}^3(x^2y),
\]
\[
v_{11}=\text{\rm tr}(x^3)\text{\rm tr}(x^2)\text{\rm tr}(xy)\text{\rm tr}(y^2)
-\text{\rm tr}(x^3)\text{\rm tr}^3(xy)
\]
\[
-\text{\rm tr}(x^2y)\text{\rm tr}^2(x^2)\text{\rm tr}(y^2)
+\text{\rm tr}(x^2y)\text{\rm tr}(x^2)\text{\rm tr}^2(xy),
\]
\[
v_{12}=-3\text{\rm tr}(x^3)\text{\rm tr}(x^2)\text{\rm tr}(xy)\text{\rm tr}(y^2)
+2\text{\rm tr}(x^3)\text{\rm tr}^3(xy)
\]
\[
+3\text{\rm tr}(x^2y)\text{\rm tr}^2(x^2)\text{\rm tr}(y^2)
-3\text{\rm tr}(xy^2)\text{\rm tr}^2(x^2)\text{\rm tr}(xy)
+\text{\rm tr}(y^3)\text{\rm tr}^3(x^2).
\]

For $\lambda=(5,4)$:
\[
w_1-v_3=0,
\]
\[
24w_2+12v_1-6v_2+12v_3-6v_4-2v_5+2v_6+v_7=0,
\]
\[
w_3+v_1=0,
\]
\[
24w_4+6v_2+12v_4-4v_6+v_7+3v_8=0,
\]
\[
v_1=\text{\rm tr}([x,y]^3x)\text{\rm tr}(xy)-\text{\rm tr}([x,y]^3y)\text{\rm tr}(x^2),
\]
\[
v_2=\text{\rm tr}([x,y]^2x^2)\text{\rm tr}(xy^2)
\]
\[
-(\text{\rm tr}([x,y]^2xy)+\text{\rm tr}([x,y]^2yx))\text{\rm tr}(x^2y)
+\text{\rm tr}([x,y]^2y^2)\text{\rm tr}(x^3),
\]
\[
v_3=\frac{1}{2}\text{\rm tr}([x,y]^2)\text{\rm tr}([x,y]^2x),
\]
\[
v_4=\text{\rm tr}([x,y]^2x)(\text{\rm tr}(x^2)\text{\rm tr}(y^2)-\text{\rm tr}^2(xy)),
\]
\[
v_5=(3\text{\rm tr}(x^4)\text{\rm tr}(xy^2)
-6\text{\rm tr}(x^3y)\text{\rm tr}(x^2y)
+(2\text{\rm tr}(x^2y^2)+\text{\rm tr}(xyxy))\text{\rm tr}(x^3))\text{\rm tr}(y^2)
\]
\[
+(-\text{\rm tr}(x^4)\text{\rm tr}(y^3)
-3\text{\rm tr}(x^3y)\text{\rm tr}(xy^2)
\]
\[
+3(2\text{\rm tr}(x^2y^2)+\text{\rm tr}(xyxy))\text{\rm tr}(x^2y)
-5\text{\rm tr}(xy^3)\text{\rm tr}(x^3))\text{\rm tr}(xy)
\]
\[
+(\text{\rm tr}(x^3y)\text{\rm tr}(y^3)
-3\text{\rm tr}(xy^3)\text{\rm tr}(x^2y)
+2\text{\rm tr}(y^4)\text{\rm tr}(x^3))\text{\rm tr}(x^2),
\]
\[
v_6=(\text{\rm tr}(x^4)\text{\rm tr}(y^3)-3\text{\rm tr}(x^3y)\text{\rm tr}(xy^2)
\]
\[
+(2\text{\rm tr}(x^2y^2)+\text{\rm tr}(xyxy))\text{\rm tr}(x^2y)
-\text{\rm tr}(xy^3)\text{\rm tr}(x^3))\text{\rm tr}(xy)
\]
\[
+(-\text{\rm tr}(x^3y)\text{\rm tr}(y^3)
+(2\text{\rm tr}(x^2y^2)+\text{\rm tr}(xyxy))\text{\rm tr}(xy^2)
\]
\[
-3\text{\rm tr}(xy^3)\text{\rm tr}(x^2y)+\text{\rm tr}(y^4)\text{\rm tr}(x^3))\text{\rm tr}(x^2),
\]
\[
v_7=\frac{1}{2}\text{\rm tr}([x,y]^2)(\text{\rm tr}(x^3)\text{\rm tr}(y^2)
-2\text{\rm tr}(x^2y)\text{\rm tr}(xy)+\text{\rm tr}(xy^2)\text{\rm tr}(x^2)),
\]
\[
v_8=(\text{\rm tr}(x^2)\text{\rm tr}(y^2)-\text{\rm tr}^2(xy))
(\text{\rm tr}(x^3)\text{\rm tr}(y^2)
-2\text{\rm tr}(x^2y)\text{\rm tr}(xy)+\text{\rm tr}(xy^2)\text{\rm tr}(x^2)).
\]
\end{lemma}

\begin{lemma}
For $\lambda=(8,2)$ the elements $w_i$ and $v_j$
satisfy the following relations in $C_0$:
\[
360w_1-90v_1-45v_2-150v_3-90v_6-4v_7-36v_8
\]
\[
-80v_9-45v_{12}+v_{13}-10v_{14}+9v_{15}=0,
\]
\[
720w_2+90v_1+90v_2-420v_3+90v_4-90v_5-135v_6+90v_8
\]
\[
-80v_9+12v_{10}+9v_{11}+60v_{12}-2v_{13}-20v_{14}-18v_{15}-45v_{16}=0,
\]
\[
720w_3-90v_1-180v_2+90v_4+150v_5-195v_6+8v_7
\]
\[
+162v_8-12v_{10}-99v_{11}-150v_{12}+40v_{14}+45v_{16}=0,
\]
\[
1440w_4+270v_2-960v_3-240v_5-120v_6-16v_7-144v_8
\]
\[
-320v_9+30v_{10}+180v_{11}+195v_{12}-80v_{14}+90v_{15}-135v_{16}=0,
\]
where
\[
v_1=\text{\rm tr}([x,y]^2x^2)\text{\rm tr}(x^4),
\]
\[
v_2=\text{\rm tr}([x,y]^2x^2)\text{\rm tr}^2(x^2),
\]
\[
v_3=\text{\rm tr}([x,y]^2x)\text{\rm tr}(x^3)\text{\rm tr}(x^2),
\]
\[
v_4=\text{\rm tr}^2(x^3y)\text{\rm tr}(x^2)
-2\text{\rm tr}(x^4)\text{\rm tr}(x^3y)\text{\rm tr}(xy)
+\text{\rm tr}^2(x^4)(y^2),
\]
\[
v_5=(2\text{\rm tr}(x^2y^2)+\text{\rm tr}(xyxy))\text{\rm tr}(x^4)\text{\rm tr}(x^2)
-3\text{\rm tr}^2(x^3y)\text{\rm tr}(x^2),
\]
\[
v_6=\frac{1}{2}\text{\rm tr}([x,y]^2)\text{\rm tr}(x^4)\text{\rm tr}(x^2),
\]
\[
v_7=5(2\text{\rm tr}(x^2y^2)+\text{\rm tr}(xyxy))\text{\rm tr}^2(x^3)
-30\text{\rm tr}(x^3y)\text{\rm tr}(x^2y)\text{\rm tr}(x^3)
\]
\[
+3\text{\rm tr}(x^4)(2\text{\rm tr}(xy^2)\text{\rm tr}(x^3)+3\text{\rm tr}^2(x^2y)),
\]
\[
v_8=\text{\rm tr}(x^4)(\text{\rm tr}(xy^2)\text{\rm tr}(x^3)-\text{\rm tr}^2(x^2y)),
\]
\[
v_9=\frac{1}{2}\text{\rm tr}([x,y]^2)\text{\rm tr}^2(x^3),
\]
\[
v_{10}=(5(2\text{\rm tr}(x^2y^2)+\text{\rm tr}(xyxy))\text{\rm tr}^2(x^2)
-30\text{\rm tr}(x^3y)\text{\rm tr}(xy)\text{\rm tr}(x^2)
\]
\[
+3\text{\rm tr}(x^4)(4\text{\rm tr}^2(xy)+\text{\rm tr}(y^2)\text{\rm tr}(x^2)))\text{\rm tr}(x^2),
\]
\[
v_{11}=\text{\rm tr}(x^4)(\text{\rm tr}(y^2)\text{\rm tr}(x^2)-\text{\rm tr}^2(xy))\text{\rm tr}(x^2),
\]
\[
v_{12}=(\text{\rm tr}(x^2y^2)-\text{\rm tr}(xyxy))\text{\rm tr}^3(x^2),
\]
\[
v_{13}=5\text{\rm tr}^2(x^3)\text{\rm tr}(x^3)(\text{\rm tr}(y^2)\text{\rm tr}(x^2)+2\text{\rm tr}^2(xy))
-30\text{\rm tr}(x^3)\text{\rm tr}(x^2y)\text{\rm tr}(xy)\text{\rm tr}(x^2)
\]
\[
+3(2\text{\rm tr}(x^3)\text{\rm tr}(xy^2)+3\text{\rm tr}^2(x^2y))\text{\rm tr}^2(x^2),
\]
\[
v_{14}=\text{\rm tr}^2(x^3)(\text{\rm tr}(x^2)\text{\rm tr}(y^2)-\text{\rm tr}^2(xy)),
\]
\[
v_{15}=(\text{\rm tr}(x^3)\text{\rm tr}(xy^2)-\text{\rm tr}^2(x^2y))\text{\rm tr}^2(x^2),
\]
\[
v_{16}=(\text{\rm tr}(x^2)\text{\rm tr}(y^2)-\text{\rm tr}^2(xy))\text{\rm tr}^3(x^2).
\]
\end{lemma}

\begin{lemma}
In $C_0$ and for $\lambda=(7,3)$, the elements $w_1,\ldots,w_7$ and
\[
v_1=\text{\rm tr}([x,y]^3x^2)\text{\rm tr}(x^2),
\]
\[
v_2=\text{\rm tr}([x,y]^3x)\text{\rm tr}(x^3),
\]
\[
v_3=(\text{\rm tr}([x,y]^2xy)+\text{\rm tr}([x,y]^2yx))\text{\rm tr}(x^4)
-2\text{\rm tr}([x,y]^2x^2)\text{\rm tr}(x^3y),
\]
\[
v_4=(\text{\rm tr}(x^2y^2xy)-\text{\rm tr}(y^2x^2yx))\text{\rm tr}(x^4),
\]
\[
v_5=((\text{\rm tr}([x,y]^2xy)+\text{\rm tr}([x,y]^2yx))\text{\rm tr}(x^2)
-2\text{\rm tr}([x,y]^2x^2)\text{\rm tr}(xy))\text{\rm tr}(x^2),
\]
\[
v_6=(\text{\rm tr}(x^2y^2xy)-\text{\rm tr}y^2x^2yx))\text{\rm tr}^2(x^2),
\]
\[
v_7=(\text{\rm tr}([x,y]^2x)\text{\rm tr}(x^2y)-\text{\rm tr}([x,y]^2y)\text{\rm tr}(x^3))\text{\rm tr}(x^2),
\]
\[
v_8=\text{\rm tr}([x,y]^2x)(\text{\rm tr}(x^3)\text{\rm tr}(xy)
-\text{\rm tr}(x^2y)\text{\rm tr}(x^2)),
\]
\[
v_9=(-6\text{\rm tr}^2(x^3y)
+2\text{\rm tr}(x^4)(2\text{\rm tr}(x^2y^2)+\text{\rm tr}(xyxy)))\text{\rm tr}(xy)
\]
\[
+(\text{\rm tr}(x^3y)(2\text{\rm tr}(x^2y^2)+\text{\rm tr}(xyxy))
-3\text{\rm tr}(x^4)\text{\rm tr}(xy^3))\text{\rm tr}(x^2),
\]
\[
v_{10}=\frac{1}{2}\text{\rm tr}([x,y]^2)(\text{\rm tr}(x^4)\text{\rm tr}(xy)-\text{\rm tr}(x^3y)\text{\rm tr}(x^2)),
\]
\[
v_{11}=\text{\rm tr}(x^4)\text{\rm tr}(x^3)\text{\rm tr}(y^3)
-3\text{\rm tr}(x^4)\text{\rm tr}(x^2y)\text{\rm tr}(xy^2)
+6\text{\rm tr}(x^3y)\text{\rm tr}^2(x^2y)
\]
\[
-2(2\text{\rm tr}(x^2y^2)+\text{\rm tr}(xyxy))\text{\rm tr}(x^3)\text{\rm tr}(x^2y)
+2\text{\rm tr}(xy^3)\text{\rm tr}^2(x^3),
\]
\[
v_{12}=-\text{\rm tr}(x^4)\text{\rm tr}(x^2y)\text{\rm tr}(xy^2)
+\text{\rm tr}(x^3y)\text{\rm tr}(x^3)\text{\rm tr}(xy^2)
+2\text{\rm tr}(x^3y)\text{\rm tr}^2(x^2y)
\]
\[
-(2\text{\rm tr}(x^2y^2)+\text{\rm tr}(xyxy))\text{\rm tr}(x^3)\text{\rm tr}(x^2y)
+\text{\rm tr}(xy^3)\text{\rm tr}^2(x^3),
\]
\[
v_{13}=-\text{\rm tr}(x^4)\text{\rm tr}(x^2)\text{\rm tr}(xy)\text{\rm tr}(y^2)
+\text{\rm tr}(x^4)\text{\rm tr}^3(xy)
\]
\[
+\text{\rm tr}(x^3y)\text{\rm tr}^2(x^2)\text{\rm tr}(y^2)
-\text{\rm tr}(x^3y)\text{\rm tr}(x^2)\text{\rm tr}^2(xy),
\]
\[
v_{14}=\text{\rm tr}(x^4)\text{\rm tr}^3(xy-3\text{\rm tr}(x^3y)\text{\rm tr}(x^2)\text{\rm tr}^2(xy)
\]
\[
+(2\text{\rm tr}(x^2y^2)+\text{\rm tr}(xyxy))\text{\rm tr}^2(x^2)\text{\rm tr}(xy)
-\text{\rm tr}(xy^3)\text{\rm tr}^3(x^2),
\]
\[
v_{15}=\text{\rm tr}^2(x^3)\text{\rm tr}(xy)\text{\rm tr}(y^2)
-\text{\rm tr}(x^3)\text{\rm tr}(x^2y)\text{\rm tr}(x^2)\text{\rm tr}(y^2)
-2\text{\rm tr}(x^3)\text{\rm tr}(x^2y)\text{\rm tr}^2(xy)
\]
\[
+\text{\rm tr}(x^3)\text{\rm tr}(xy^2)\text{\rm tr}(x^2)\text{\rm tr}(xy)
+2\text{\rm tr}^2(x^2y)\text{\rm tr}(x^2)\text{\rm tr}(xy)
-\text{\rm tr}(x^2y)\text{\rm tr}(xy^2)\text{\rm tr}^2(x^2),
\]
\[
v_{16}=-2\text{\rm tr}(x^3)\text{\rm tr}(xy^2)\text{\rm tr}(x^2)\text{\rm tr}(xy)
+2\text{\rm tr}^2(x^2y)\text{\rm tr}(x^2)\text{\rm tr}(xy)
\]
\[
+\text{\rm tr}(x^3)\text{\rm tr}(y^3)\text{\rm tr}^2(x^2-\text{\rm tr}(x^2y)\text{\rm tr}(xy^2)\text{\rm tr}^2(x^2),
\]
satisfy the relations
\[
24w_1-12v_1-8v_2-18v_4+9v_6=0,
\]
\[
48w_2-24v_2+12v_3-24v_4-3v_5-8v_7+4v_8-12v_{10}+4v_{11}-2v_{15}-v_{16}=0,
\]
\[
48w_3-24v_1-8v_2+12v_3-12v_4-3v_5+6v_6-8v_7
\]
\[
+4v_8-12v_{10}+4v_{11}-2v_{15}-v_{16}=0,
\]
\[
12w_4-4v_2-3v_6=0,
\]
\[
48w_5-60v_1+8v_2+6v_3-24v_4+36v_6+10v_7
\]
\[
+4v_8-4v_9-10v_{10}+10v_{11}-12v_{12}+6v_{13}-2v_{15}-v_{16}=0,
\]
\[
48w_6-48v_1-8v_2-36v_4-3v_5+30v_6-32v_7
\]
\[
-20v_8-8v_{11}+24v_{12}-6v_{13}+6v_{14}-6v_{15}+3v_{16}=0,
\]
\[
12w_7-6v_1+4v_2+3v_6=0.
\]
\end{lemma}

\begin{lemma}
The relation
\[
\sum_{i=1}^{10}\alpha_iw_i+\sum_{j=1}^{24}\beta_jv_j=0,
\]
holds in $C_0$ between $w_i$ and $v_j$ for $\lambda=(6,4)$, where
\[
v_1=(\text{\rm tr}([x,y]^3xy)+\text{\rm tr}([x,y]^3yx))\text{\rm tr}(x^2)
-2\text{\rm tr}([x,y]^3x^2)\text{\rm tr}(xy),
\]
\[
v_2=\text{\rm tr}([x,y]^4)\text{\rm tr}(x^2),
\]
\[
v_3=\text{\rm tr}([x,y]^3x)\text{\rm tr}(x^2y)-\text{\rm tr}([x,y]^3y)\text{\rm tr}(x^3),
\]
\[
v_4=\text{\rm tr}([x,y]^2x^2)(2\text{\rm tr}(x^2y^2)+\text{\rm tr}(xyxy))
\]
\[
-3(\text{\rm tr}([x,y]^2xy)+\text{\rm tr}([x,y]^2yx))\text{\rm tr}(x^3y)
+3\text{\rm tr}([x,y]^2y^2)\text{\rm tr}(x^4),
\]
\[
v_5=\frac{1}{2}\text{\rm tr}([x,y]^2x^2)\text{\rm tr}([x,y]^2),
\]
\[
v_6=\text{\rm tr}([x,y]^2x^2)\text{\rm tr}^2(xy)
\]
\[
-(\text{\rm tr}([x,y]^2yx)+\text{\rm tr}([x,y]^2xy))\text{\rm tr}(x^2)\text{\rm tr}(xy)
+\text{\rm tr}([x,y]^2y^2)\text{\rm tr}^2(x^2),
\]
\[
v_7=\text{\rm tr}([x,y]^2x^2)(\text{\rm tr}(x^2)\text{\rm tr}(y^2)-\text{\rm tr}^2(xy)),
\]
\[
v_8=\text{\rm tr}^2t([x,y]^2x),
\]
\[
v_9=\text{\rm tr}([x,y]^2x)(\text{\rm tr}(x^3)\text{\rm tr}(y^2)
-2\text{\rm tr}(x^2y)\text{\rm tr}(xy)+\text{\rm tr}(xy^2)\text{\rm tr}(x^2)),
\]
\[
v_{10}=\text{\rm tr}([x,y]^2x)(-\text{\rm tr}(x^3)\text{\rm tr}(y^2)+\text{\rm tr}(x^2y)\text{\rm tr}(xy))
\]
\[
+\text{\rm tr}([x,y]^2y)(\text{\rm tr}(x^3)\text{\rm tr}(xy)-\text{\rm tr}(x^2y)\text{\rm tr}(x^2)),
\]
\[
v_{11}=(\text{\rm tr}(x^4)\text{\rm tr}(x^2y^2)+\text{\rm tr}(x^4)\text{\rm tr}(xyxy)
-2\text{\rm tr}^2(x^3y))\text{\rm tr}(y^2)
\]
\[
+(-2\text{\rm tr}(x^4)\text{\rm tr}(xy^3)
+2\text{\rm tr}(x^3y)\text{\rm tr}(x^2y^2))\text{\rm tr}(xy)
\]
\[
+(\text{\rm tr}(x^4)\text{\rm tr}(y^4)
-2\text{\rm tr}(x^3y)\text{\rm tr}(xy^3)+\text{\rm tr}(x^2y^2)\text{\rm tr}(xyxy))\text{\rm tr}(x^2),
\]
\[
v_{12}=(\text{\rm tr}(x^4)\text{\rm tr}(x^2y^2)-\text{\rm tr}^2(x^3y))\text{\rm tr}(y^2)
\]
\[
+(-\text{\rm tr}(x^4)\text{\rm tr}(xy^3)+\text{\rm tr}(x^3y)\text{\rm tr}(xyxy))\text{\rm tr}(xy)
\]
\[
+(\text{\rm tr}(x^3y)\text{\rm tr}(xy^3)-\text{\rm tr}(x^2y^2)\text{\rm tr}(xyxy))\text{\rm tr}(x^2),
\]
\[
v_{13}=(\text{\rm tr}(x^4)\text{\rm tr}(y^4)
-4\text{\rm tr}(x^3y)\text{\rm tr}(xy^3)+\text{\rm tr}^2(x^2y^2)
+2\text{\rm tr}(x^2y^2)\text{\rm tr}(xyxy))\text{\rm tr}(x^2),
\]
\[
v_{14}=(-\text{\rm tr}(x^4)\text{\rm tr}(y^4)+4\text{\rm tr}(x^3y)\text{\rm tr}(xy^3)
-4\text{\rm tr}(x^2y^2)\text{\rm tr}(xyxy)+\text{\rm tr}^2(xyxy))\text{\rm tr}(x^2),
\]
\[
v_{15}=\text{\rm tr}(x^4\text{\rm tr}(x^2y)\text{\rm tr}(y^3)
+\text{\rm tr}(x^3y)(-\text{\rm tr}(x^3)\text{\rm tr}(y^3)-3\text{\rm tr}(x^2y)\text{\rm tr}(xy^2))
\]
\[
+\text{\rm tr}(x^2y^2)(3\text{\rm tr}(x^3)\text{\rm tr}(xy^2)+\text{\rm tr}^2(x^2y))
\]
\[
+2\text{\rm tr}(xyxy)\text{\rm tr}^2(x^2y)-4\text{\rm tr}(xy^3)\text{\rm tr}(x^3)\text{\rm tr}(x^2y)
+\text{\rm tr}(y^4)\text{\rm tr}^2(x^3),
\]
\[
v_{16}=\text{\rm tr}(x^4)\text{\rm tr}^2(xy^2)-4\text{\rm tr}(x^3y)\text{\rm tr}(x^2y)\text{\rm tr}(xy^2)
+\text{\rm tr}(x^2y^2)(2\text{\rm tr}(x^3)\text{\rm tr}(xy^2)
\]
\[
+2\text{\rm tr}^2(x^2y))
+2\text{\rm tr}(xyxy)\text{\rm tr}^2(x^2y)-4\text{\rm tr}(xy^3)\text{\rm tr}(x^3)\text{\rm tr}(x^2y)
+\text{\rm tr}(y^4)\text{\rm tr}^2(x^3),
\]
\[
v_{17}=\text{\rm tr}(x^2y^2)(-\text{\rm tr}(x^3)\text{\rm tr}(xy^2)+\text{\rm tr}^2(x^2y))
+\text{\rm tr}(xyxy)(\text{\rm tr}(x^3)\text{\rm tr}(xy^2)-\text{\rm tr}^2(x^2y)),
\]
\[
v_{18}=\text{\rm tr}(x^4)(\text{\rm tr}(x^2)\text{\rm tr}(y^2)-\text{\rm tr}^2(xy))\text{\rm tr}(y^2)
\]
\[
+\text{\rm tr}(x^3y)(-2\text{\rm tr}(x^2)\text{\rm tr}(y^2)+2\text{\rm tr}^2(xy))\text{\rm tr}(xy)
\]
\[
+\text{\rm tr}(x^2y^2)\text{\rm tr}(x^2)(\text{\rm tr}(x^2)\text{\rm tr}(y^2)-\text{\rm tr}^2(xy)),
\]
\[
v_{19}=\text{\rm tr}(x^4)(-\text{\rm tr}(x^2)\text{\rm tr}(y^2)+\text{\rm tr}^2(xy))\text{\rm tr}(y^2)
\]
\[
+\text{\rm tr}(x^3y)(2\text{\rm tr}(x^2)\text{\rm tr}(y^2)-2\text{\rm tr}^2(xy))\text{\rm tr}(xy)
\]
\[
+\text{\rm tr}(xyxy)\text{\rm tr}(x^2)(-\text{\rm tr}(x^2)\text{\rm tr}(y^2)+\text{\rm tr}^2(xy)),
\]
\[
v_{20}=\text{\rm tr}(x^4)\text{\rm tr}(x^2\text{\rm tr}^2(y^2)
-4\text{\rm tr}(x^3y)\text{\rm tr}(x^2)\text{\rm tr}(xy)\text{\rm tr}(y^2)
+4\text{\rm tr}(x^2y^2)\text{\rm tr}(x^2)\text{\rm tr}^2(xy)
\]
\[
+2\text{\rm tr}(xyxy)\text{\rm tr}^2(x^2)\text{\rm tr}(y^2)
-4\text{\rm tr}(xy^3)\text{\rm tr}^2(x^2)\text{\rm tr}(xy)
+\text{\rm tr}(y^4)\text{\rm tr}^3(x^2),
\]
\[
v_{21}=\text{\rm tr}^2(x^3)\text{\rm tr}^2(y^2)
-4\text{\rm tr}(x^3)\text{\rm tr}(x^2y)\text{\rm tr}(xy)\text{\rm tr}(y^2)
\]
\[
+2\text{\rm tr}^2(x^2y)\text{\rm tr}(x^2)\text{\rm tr}(y^2)
+4\text{\rm tr}(x^3)\text{\rm tr}(xy^2)\text{\rm tr}^2(xy)
+(-2\text{\rm tr}(x^3)\text{\rm tr}(y^3)
\]
\[
-2\text{\rm tr}(x^2y)\text{\rm tr}(xy^2))\text{\rm tr}(x^2)\text{\rm tr}(xy)
+(2\text{\rm tr}(x^2y)\text{\rm tr}(y^3)-\text{\rm tr}^2(xy^2))\text{\rm tr}^2(x^2),
\]
\[
v_{22}=((\text{\rm tr}(x^3)\text{\rm tr}(xy^2)-\text{\rm tr}^2(x^2y))\text{\rm tr}(y^2)
\]
\[
+(-\text{\rm tr}(x^3)\text{\rm tr}(y^3)+\text{\rm tr}(x^2y)\text{\rm tr}(xy^2))\text{\rm tr}(xy)
\]
\[
+(\text{\rm tr}(x^2y)\text{\rm tr}(y^3)-\text{\rm tr}^2(xy^2)))\text{\rm tr}(x^2),
\]
\[
v_{23}=(-\text{\rm tr}(x^3)\text{\rm tr}(xy^2)+\text{\rm tr}^2(x^2y))\text{\rm tr}^2(xy)
\]
\[
+(\text{\rm tr}(x^3)\text{\rm tr}(y^3)-\text{\rm tr}(x^2y)\text{\rm tr}(xy^2))\text{\rm tr}(x^2)\text{\rm tr}(xy)
+(-\text{\rm tr}(x^2y)\text{\rm tr}(y^3)+\text{\rm tr}^2(xy^2))\text{\rm tr}^2(x^2),
\]
\[
v_{24}=(\text{\rm tr}(x^2)\text{\rm tr}(y^2)-\text{\rm tr}^2(xy))^2\text{\rm tr}(x^2),
\]
and
\[
\beta_1=\frac{1}{4}(\alpha_3+\alpha_6-\alpha_7-\alpha_9),
\]
\[
\beta_2=\frac{1}{24}(-6\alpha_1-3\alpha_2-8\alpha_3-8\alpha_4-17\alpha_5+2\alpha_6-5\alpha_7
-\alpha_8-2\alpha_9-6\alpha_{10}),
\]
\[
\beta_3=\frac{1}{2}(\alpha_3+\alpha_4-\alpha_6+\alpha_7-2\alpha_9),
\]
\[
\beta_4=\frac{1}{12}(-\alpha_2+\alpha_3+\alpha_7-\alpha_9),
\]
\[
\beta_5=\frac{1}{12}(-12\alpha_1-2\alpha_2+5\alpha_3-3\alpha_4+3\alpha_6+2\alpha_7+4\alpha_9-3\alpha_{10}),
\]
\[
\beta_6=\frac{1}{48}(-3\alpha_2-28\alpha_3-16\alpha_4-37\alpha_5+4\alpha_6-13\alpha_7
-5\alpha_8-10\alpha_9-12\alpha_{10}),
\]
\[
\beta_7=\frac{1}{48}(-15\alpha_2-40\alpha_3-16\alpha_4-37\alpha_5+4\alpha_6
-13\alpha_7-5\alpha_8-22\alpha_9-24\alpha_{10}),
\]
\[
\beta_8=\frac{1}{4}(-\alpha_2+\alpha_3+3\alpha_4+3\alpha_5-\alpha_6-\alpha_7+\alpha_8-3\alpha_{10}),
\]
\[
\beta_9=\frac{1}{12}(+3\alpha_2+5\alpha_3+2\alpha_4+5\alpha_5+\alpha_6+2\alpha_7+\alpha_8-\alpha_9),
\]
\[
\beta_{10}=\frac{1}{12}(-2\alpha_3+\alpha_4-2\alpha_5+2\alpha_6-2\alpha_7-4\alpha_8+\alpha_9-6\alpha_{10}),
\]
\[
\beta_{11}=\frac{1}{8}(+\alpha_2-\alpha_3-2\alpha_4+\alpha_7-\alpha_9+2\alpha_{10}),
\]
\[
\beta_{12}=\frac{1}{4}(-\alpha_2+7\alpha_3+2\alpha_4
-2\alpha_6+\alpha_7+3\alpha_9+2\alpha_{10}),
\]
\[
\beta_{13}=\frac{1}{24}(+12\alpha_1-6\alpha_2+13\alpha_3+7\alpha_4-14\alpha_5-7\alpha_6
-2\alpha_7-10\alpha_8+10\alpha_9+3\alpha_{10}),
\]
\[
\beta_{14}=\frac{1}{24}(+12\alpha_1+6\alpha_2+8\alpha_3
+11\alpha_4+17\alpha_5-5\alpha_6+2\alpha_7
+\alpha_8+5\alpha_9+9\alpha_{10}),
\]
\[
\beta_{15}=\frac{1}{12}(-\alpha_2-\alpha_3-2\alpha_4-8\alpha_5+2\alpha_6-\alpha_7-4\alpha_8+3\alpha_9),
\]
\[
\beta_{16}=\frac{1}{4}(+\alpha_4+3\alpha_5-\alpha_6+\alpha_8-\alpha_9-\alpha_{10}),
\]
\[
\beta_{17}=\frac{1}{12}(+\alpha_2-2\alpha_3-\alpha_4+2\alpha_5
+\alpha_6+\alpha_7-2\alpha_8-3\alpha_9-3\alpha_{10}),
\]
\[
\beta_{18}=\frac{1}{48}(-15\alpha_2-64\alpha_3-22\alpha_4-43\alpha_5
+10\alpha_6-19\alpha_7-11\alpha_8-28\alpha_9-36\alpha_{10}),
\]
\[
\beta_{19}=\frac{1}{48}(-9\alpha_2-46\alpha_3-22\alpha_4
-43\alpha_5+10\alpha_6-13\alpha_7-11\alpha_8-22\alpha_9-24\alpha_{10}),
\]
\[
\beta_{20}=\frac{1}{96}(+9\alpha_2+4\alpha_3+10\alpha_4+25\alpha_5
+2\alpha_6+\alpha_7+5\alpha_8+4\alpha_9),
\]
\[
\beta_{21}=\frac{1}{48}(+3\alpha_2+9\alpha_3-\alpha_4+3\alpha_5
-\alpha_6+5\alpha_7+5\alpha_8+11\alpha_{10}),
\]
\[
\beta_{22}=\frac{1}{48}(-3\alpha_2+18\alpha_3+4\alpha_4+15\alpha_5
-14\alpha_6+7\alpha_7+7\alpha_8+10\alpha_{10}),
\]
\[
\beta_{23}=\frac{1}{4}(+2\alpha_3+\alpha_5
-\alpha_6+\alpha_7+\alpha_8+2\alpha_{10}),
\]
\[
\beta_{24}=\frac{1}{16}(+\alpha_2+2\alpha_3+\alpha_5
+\alpha_7+\alpha_8+2\alpha_{10}).
\]
\end{lemma}

\begin{lemma}
For $\lambda=(5,5)$ in $C_0$:
\[
w_1-5v_2=0,
\]
\[
w_2-w_3+v_1-4v_2-v_3=0,
\]
\[
w_3+w_4-\frac{5}{2}v_1-v_2+2v_3=0,
\]
where
\[
v_1=\text{\rm tr}([x,y]^3x^2)\text{\rm tr}(y^2)
-(\text{\rm tr}([x,y]^3yx)+\text{\rm tr}([x,y]^3xy)\text{\rm tr}(xy)
+\text{\rm tr}([x,y]^3y^2)\text{\rm tr}(x^2),
\]
\[
v_2=\frac{1}{2}(\text{\rm tr}(x^2y^2xy)-\text{\rm tr}(y^2x^2yx))
\text{\rm tr}([x,y]^2),
\]
\[
v_3=(\text{\rm tr}(x^2y^2xy)-\text{\rm tr}(y^2x^2yx))
(\text{\rm tr}(x^2)\text{\rm tr}(y^2)-\text{\rm tr}^2(xy)),
\]
\[
v_4=(\text{\rm tr}(x^3)\text{\rm tr}(y^2)
-2\text{\rm tr}(x^2y)\text{\rm tr}(xy)
+\text{\rm tr}(xy^2)\text{\rm tr}(x^2))
\text{\rm tr}([x,y]^2y)
\]
\[
-(\text{\rm tr}(y^3)\text{\rm tr}(x^2)
-2\text{\rm tr}(xy^2)\text{\rm tr}(xy)
+\text{\rm tr}(x^2y)\text{\rm tr}(y^2))
\text{\rm tr}([x,y]^2x).
\]
\end{lemma}

The following theorem is the main result of the paper.

\begin{theorem} The minimal generating $GL_2$-module $G$ of the pure trace
algebra $C_{42}$ of two generic $4\times 4$ matrices decomposes as
\[
\begin{array}{c}
G=W(1,0)\oplus W(2,0)\oplus W(3,0)\oplus W(4,0)\oplus W(2,2)\oplus W(3,2)\\
\oplus W(4,2)\oplus W(3,3)\oplus W(4,3)\oplus W(5,3)\oplus W(4,4)\oplus W(6,3)
\oplus W(5,5).\\
\end{array}
\]
For all
$\lambda=(\lambda_1,\lambda_2)\not=(5,5)$, the module
$W(\lambda_1,\lambda_2)$ is generated by the element
\[
w_{\lambda}(X,Y)=\text{\rm tr}((XY-YX)^{\lambda_2}X^{\lambda_1-\lambda_2}).
\]
For $\lambda=(5,5)$ a generator of $W(5,5)$ is
\[
w_{(5,5)}(X,Y)=\text{\rm tr}((XY-YX)^3(X^2Y^2-XYYX-YXXY+Y^2X^2)).
\]
\end{theorem}

\begin{proof}
Using (\ref{from generic to traceless matrices}), we may replace
the elements $\text{\rm tr}(Z_1\cdots Z_k)$, $Z_i=X,Y$, $k\geq 2$,
modulo the subalgebra generated by products of traces of degree $<k$,
with $\text{\rm tr}(z_1\cdots z_k)$, where $z_i=x,y$ are
generic traceless matrices. Hence, it is
sufficient to show that the algebra $C_0$ has a minimal $GL_2$-module
of generators
\begin{equation}
\begin{array}{c}\label{generators of C0}
M=W(2,0)\oplus W(3,0)\oplus W(4,0)\oplus W(2,2)\oplus W(3,2)\\
\oplus W(4,2)\oplus W(3,3)\oplus W(4,3)\oplus W(5,3)\oplus W(4,4)\oplus W(6,3)
\oplus W(5,5),\\
\end{array}
\end{equation}
each $W(\lambda_1,\lambda_2)$,
$(\lambda_1,\lambda_2)\not=(5,5)$, is generated by
\[
w_{\lambda}(x,y)=\text{\rm tr}((xy-yx)^{\lambda_2}x^{\lambda_1-\lambda_2}),
\]
and for $\lambda=(5,5)$ a generator of $W(5,5)$ is
\begin{equation}\label{generator of degree 10}
w_{(5,5)}(x,y)=\text{\rm tr}((xy-yx)^3(x^2y^2-xyyx-yxxy+y^2x^2)).
\end{equation}
There are three typical cases. We consider the image of
the $GL_2$-modules $U_n$ in $C_0$ under the natural homomorphism
$\pi:P(x,y)\to C_0$. We may assume that the generators of $C_0$
are in $\pi(U_2\oplus\cdots\oplus U_{10})$. Lemma 3.1 allows
to consider the cases $n\geq 6$ only.

Let $\lambda=(n)$, $n=6,\ldots,10$. This case was handled
in Lemma 3.1. The only $GL_2$-module
$W(n)$ in $\pi(U_n)$ is generated by $\text{tr}(x^n)$.
By (\ref{trace of xxxxx}), $\text{tr}(x^5)$ belongs to the subalgebra
of $C_0$ generated by traces of lower degree.
Similar equations hold for the other $n\geq 6$. Hence, $W(n)$ does not
participate in any minimal $GL_2$-module of generators of $C_0$
for $n=6,\ldots,10$.

Let $\lambda=(6,2)$. By Proposition 2.1, $W(6,2)$ participates
with multiplicity 3 in $U_8$. The highest weight vectors of the
three copies of $W(6,2)$ are $w_1,w_2,w_3$ from Lemma 2.2. By
the relations in Lemma 3.4,
$\pi(w_1),\pi(w_2),\pi(w_3)$ are linear combinations of
$v_1,\ldots,v_{10}$, and hence belong to the subalgebra of $C_0$
generated by the elements of degree $\leq 7$
(and even of degree $\leq 6$). Therefore, as in the case $\lambda=(n)$,
the minimal $GL_2$-modules of generators of $C_0$ do not contain
$W(6,2)$. The cases $\lambda=(5,2), (7,2), (5,4), (8,2), (7,3), (6,4)$
are similar.

Let $\lambda=(5,3)$. By Proposition 2.1 again, $W(5,3)$ participates
with multiplicity 3 in $U_8$. The corresponding highest weight vectors
are $w_1,w_2,w_3$ from Lemma 2.2, and the element $w_1$ is equal to
$\text{\rm tr}((xy-yx)^3x^2)$. By the relations in Lemma 3.4,
$\pi(w_2),\pi(w_3)$ are linear combinations of $\pi(w_1)$ and
$v_1,\ldots,v_5$. Hence they belong to the subalgebra of $C_0$
generated by the elements of degree $\leq 7$ and by the $GL_2$-module
generated by $w_1$. In this way, the multiplicity of $W(5,3)$
in the minimal $GL_2$-module of generators is bounded from above by 1.
The cases $\lambda=(4,2), (4,3), (4,4), (6,3)$ are similar.
The case $\lambda=(5,5)$ is analogous but one has to express
the elements $w_1,w_3,w_4$ in terms of $w_2$ and $v_1,v_2,v_3,v_4$.
It is easy to see that $w_2$ is equal to
(\ref{generator of degree 10}).
Finally, the multiplicity of $W(3,3)$ in $U_6$ is already 1.
In this way we obtain that the minimal $GL_2$-module of generators
of $C_0$ is a homomorphism image of the module $M$ from
(\ref{generators of C0}), and the algebra $C_0$ is a homomorphic
image of the symmetric algebra $K[M]$ of $M$.

Clearly, the coefficients of the Hilbert series of $K[M]$
are bigger or equal to the coefficients of the Hilbert series of
$C_0$. The Hilbert series of $C_{42}$ is given in
(\ref{Hilbert series of C42}). The factor $(1-t)(1-u)$
in its denominator corresponds to the polynomial algebra
$K[\text{tr}(X),\text{tr}(Y)]$ in
(\ref{the tensor product}). Hence, the Hilbert series of $C_0$
is equal to $h(t,u)=P_C(t,u)/Q_C(t,u)$ from (\ref{Hilbert series of C0}).
It follows from (\ref{generators of C0}) that the Hilbert series of $K[M]$ is
\[
H(K[M],t,u)=\frac{1}{q_2(t,u)\cdots q_{10}(t,u)},
\]
where $q_2(t,u),\ldots,q_5(t,u)$ are given in
(\ref{first four polynomials}) and
\[
q_6(t,u)=(1-t^4u^2)(1-t^3u^3)^2(1-t^2u^4),
\]
\[
q_7(t,u)=(1-t^4u^3)(1-t^3u^4),
\]
\[
q_8(t,u)=(1-t^5u^3)(1-t^4u^4)^2(1-t^3u^5),
\]
\[
q_9(t,u)=(1-t^6u^3)(1-t^5u^4)(1-t^4u^5)(1-t^3u^6),\quad
q_{10}(t,u)=(1-t^5u^5).
\]
By direct calculations we see that the coefficients of the homogeneous
components of degree $\leq 10$ of $H(K[M],t,u)$ are the same as
the coefficients $h_0,h_1,\ldots,h_{10}$ from
(\ref{the first eleven coefficients}). This implies that
the generators (\ref{generators of C0}) satisfy
no relations of degree $\leq 10$ and competes the proof of the theorem.
\end{proof}

\begin{remark}
Since $\text{\rm tr}([X,Y])=0$, the equation (\ref{trace of xxxxx})
gives that the element $\text{\rm tr}([X,Y]^5)$ of $C_{42}$ belongs
to the subalgebra generated by traces of lower degree.
\end{remark}

\begin{remark}
Calculating the first homogeneous components of the difference
$H(C_0,t,u)-H(K[M],t,u)$, we see that
\[
H(C_0,t,u)-H(K[M],t,u)=
(S_{(7,5)}(t,u)+2S_{(6,6)}(t,u))
+(S_{(8,5)}(t,u)+2S_{(7,6)}(t,u))
\]
\[
+(2S_{(9,5)}(t,u)+6S_{(8,6)}(t,u)+2S_{(7,7)}(t,u))
\]
\[
+(2S_{(10,5)}(t,u)+9S_{(9,6)}(t,u)+7S_{(8,7)}(t,u)).
\]
Hence, the $GL_2$-modules of the defining relations of degree
12 and 13 are, respectively, $W(7,5)\oplus 2W(6,6)$ and
$W(8,5)\oplus 2W(7,6)$. In principle, the defining relations can
be found as those in Lemmas 3.2 -- 3.9, but the calculations
should be more complicated.
\end{remark}

\begin{remark}
The homogeneous system of parameters of $C_{42}$ found by
Teranishi \cite{T1, T2} contains all traces of degree $\leq 4$
and two elements of degree $(4,2)$ and $(2,4)$, respectively.
Translated in our language, we can choose for
a homogeneous system of parameters of $C_{42}$ the elements
\[
\text{tr}(X),\text{tr}(Y),\text{tr}(X^2),\text{tr}(XY),\text{tr}(Y^2),
\]
\[
\text{tr}(X^3),\text{tr}(X^2Y),\text{tr}(XY^2),\text{tr}(Y^3),
\]
\[
\text{tr}(X^4),\text{tr}(X^3Y),\text{tr}(X^2Y^2),\text{tr}(XYXY),
\text{tr}(XY^3),\text{tr}(Y^4),
\]
and the polynomials
\[
w_{(4,2)}(X,Y)=\text{\rm tr}((XY-YX)^2X^2),\quad
w_{(4,2)}(Y,X)=\text{\rm tr}((XY-YX)^2Y^2)
\]
from Theorem 3.10.
\end{remark}

\end{document}